\documentclass[aop,MSNbibl,dvips]{arximspdf}

\doi{10.1214/12-AOP816} %kopijuoti is PTS
\volume{41}
\issue{6}
\pubyear{2013}
\firstpage{4080}
\lastpage{4115}

\makeatletter

\newtheorem{theorem}{Theorem}
\newproclaim{remark}{Remark}
\newtheorem{lemma}{Lemma}
\newtheorem{conj}{Conjecture}
\newtheorem{prop}{Proposition}
\newtheorem{coro}{Corollary}

\newcommand{\Pois}{\operatorname{Pois}}
\newcommand{\Bi}{\operatorname{Bi}}
\newcommand{\pr}{\mathbb{P}}
\newcommand{\E}{\mathbb{E}}
\newcommand{\R}{\mathbb{R}}
\newcommand{\Z}{\mathbb{Z}}
\newcommand{\G}{\mathbb{G}}
\newcommand{\distr}{\stackrel{d}{=}}
\newcommand{\rM}{\mathcal{M}}

\makeatother

\begin{document}
\begin{frontmatter}

\title{Combinatorial approach to the interpolation method and scaling
limits in sparse random graphs}
\runtitle{Combinatorial approach to the interpolation method}

\begin{aug}
\author[A]{\fnms{Mohsen} \snm{Bayati}\thanksref{t1}\ead[label=e1]{bayati@stanford.edu}},
\author[B]{\fnms{David} \snm{Gamarnik}\corref{}\thanksref{t2}\ead[label=e2]{gamarnik@mit.edu}}
\and
\author[C]{\fnms{Prasad} \snm{Tetali}\thanksref{t3}\ead[label=e3]{tetali@math.gatech.edu}}
\runauthor{M. Bayati, D. Gamarnik and P. Tetali}
\affiliation{Stanford Univeristy, MIT and Georgia Tech}
\address[A]{M. Bayati\\
Graduate School of Business\\
Stanford University\\
655 Knight Way\\
Stanford, California 94305\\
USA\\
\printead{e1}}
\address[B]{D. Gamarnik\\
MIT Sloan School of Management\\
100 Main Street\\
Cambridge, Massachusetts 02139\\
USA\\
\printead{e2}}
\address[C]{P. Tetali\\
School of Mathematics\\
\quad and School of Computer Science\\
Georgia Institute of Technology\\
Atlanta, Georgia 30332\\
USA\\
\printead{e3}} %adresu isvedimo komanda gale!
\end{aug}

\thankstext{t1}{Supported in part by MBA Class of 1969 Faculty Scholarship.}
\thankstext{t2}{Supported in part by the NSF Grants CMMI-1031332.}
\thankstext{t3}{Supported in part by the NSF Grants DMS-07-01043 and
CCR-0910584.}

% HISTORY:
\received{\smonth{10} \syear{2011}}
\revised{\smonth{10} \syear{2012}}

% ABSTRACT
%
\begin{abstract}
We establish the existence of free energy limits for several
combinatorial models on Erd{\"o}s--R\'{e}nyi graph $\G(N,\lfloor
cN\rfloor)$ and random $r$-regular graph $\G(N,r)$. For a variety of
models, including independent sets, MAX-CUT, coloring and K-SAT, we
prove that the free energy both at a positive and zero temperature,
appropriately rescaled, converges to a limit as the size of the
underlying graph diverges to infinity. In the zero temperature case,
this is interpreted as the existence of the scaling limit for the
corresponding combinatorial optimization problem. For example, as a
special case we prove that the size of a largest independent set in
these graphs, normalized by the number of nodes converges to a limit
w.h.p. This resolves an open problem which was proposed by Aldous
(Some open problems) as one of his six favorite open problems.
It was also mentioned as an open problem in several other places:
Conjecture 2.20 in Wormald [In \textit{Surveys in Combinatorics}, 1999
(\textit{Canterbury}) (1999) 239--298 Cambridge Univ. Press];
Bollob{\'a}s and Riordan [\textit{Random Structures Algorithms}
\textbf{39} (2011) 1--38]; Janson and Thomason [\textit{Combin. Probab.
Comput.} \textbf{17} (2008) 259--264] and Aldous and Steele
[In \textit{Probability on Discrete Structures} (2004) 1--72 Springer].

Our approach is based on extending and simplifying the interpolation
method of Guerra and Toninelli [\textit{Comm. Math. Phys.} \textbf{230}
(2002) 71--79] and Franz and Leone [\textit{J. Stat. Phys.}
\textbf{111} (2003) 535--564]. Among other applications, this method
was used to prove the existence of free energy limits for Viana--Bray
and K-SAT models on Erd{\"o}s--R\'{e}nyi graphs. The case of zero
temperature was treated by taking limits of positive temperature
models. We provide instead a simpler combinatorial approach and work
with the zero temperature case (optimization) directly both in the case
of Erd{\"o}s--R\'{e}nyi graph $\G(N,\lfloor cN\rfloor)$ and random
regular graph $\G(N,r)$. In addition we establish the large deviations
principle for the satisfiability property of the constraint
satisfaction problems, coloring, K-SAT and NAE-K-SAT, for the
$\G(N,\lfloor cN\rfloor)$ random graph model.
\end{abstract}

% KEYWORDS

\end{frontmatter}

%s1 #&#
\section{Introduction}
Consider two random graph models on nodes $[N]\triangleq\{1,\ldots,\allowbreak N\}
$, the Erd{\"o}s--R\'{e}nyi graph $\G(N,M)$
and the random $r$-regular graph $\G(N,r)$.
The first model is obtained by generating $M$ edges of the $N(N-1)/2$
possible edges uniformly at random without replacement.
Specifically, assume $M=\lfloor cN\rfloor$ where $c>0$
is a constant (does not grow with $N$). The second model $\G(N,r)$ is
a graph chosen uniformly at random from the space of all $r$-regular graphs
on $N$ nodes, where the integer $r$ is a fixed integer constant.
Consider the size $|\mathcal{I}_N|$ of a largest independent set
$\mathcal{I}_N\subset[N]$ in $\G(N,\lfloor cN\rfloor)$ or $\G(N,r)$.
It is straightforward to see that $|\mathcal{I}_N|$ grows linearly
with $N$. It was conjectured in several papers including
Conjecture 2.20 in \cite{WormaldModelsRandomGraphs,gamarnikMaxWeightIndSet,BollobasRiordanMetrics}, as well as \cite{JansonThomason} and
\cite{AldousSteelesurvey}, that $|\mathcal{I}_N|/N$ converges in
probability as $N\rightarrow\infty$.
Additionally, this problem was
listed by Aldous as one of his six favorite open problems \cite
{AldousFavoriteProblems}. (For a new collection of Aldous's favorite
open problems,
see \cite{AldousFavoriteProblemsNew}.)
The fact that the actual value of $|\mathcal{I}_N|$ concentrates
around its mean
follows from a standard Azuma-type inequality. However, a real
challenge is to show that the expected value of $|\mathcal{I}_N|$
normalized by
$N$ does not fluctuate for large $N$.

This conjecture is in fact just one of a family of similar conjectures.
Consider, for example,
the random MAX-K-SAT problem---the problem of finding the largest
number of satisfiable clauses of size $K$ in a uniformly random instance
of a K-SAT problem on $N$ variables with $cN$ clauses.
This problem can be viewed as an optimization problem over a sparse random
\textit{hypergraph}. A straightforward argument shows that asymptotically
as $N\rightarrow\infty$,
at least $1-2^{-K}$ fraction of the clauses can be satisfied with high
probability (w.h.p.). Indeed any random assignment of variables
satisfies each clause
with probability $1-2^{-K}$.
It was conjectured in \cite{CopGamMohSor} that the proportion of the
largest number of satisfiable clauses has a limit w.h.p. as
$N\rightarrow\infty$.
As another example, consider the problem of partial $q$-coloring of a
graph: finding a $q$-coloring of nodes which maximizes the total
number of properly colored edges. It is natural to conjecture again
that value of this maximum has a scaling limit w.h.p. (though we are
not aware of any
papers explicitly stating this conjecture).

Recently a powerful rigorous statistical physics method was introduced
by Guerra and Toninelli \cite{GuerraTon} and further developed by
Franz and Leone \cite{FranzLeone}, Franz, Leone and Toninelli
\cite{FranzLeoneToninelliRegular}, Panchenko and Talagrand
\cite{PanchenkoTalagrand} and Montanari \cite{MontanariLDPCInterpolation}
in the context of the theory of spin glasses.
The method is based on an ingenious
interpolation between a random hypergraph model on $N$ nodes on the one
hand, and a disjoint union of random hypergraph models on $N_1$ and
$N_2$ nodes, on the other hand,
where $N=N_1+N_2$. Using this method it is possible to show for certain
spin glass models on random hypergraphs, that when one considers the
expected log-partition function,
the derivative of the interpolation function has a definite sign at
\textit{every value} of the
interpolation parameter. As a result the expected log-partition
function of the $N$-node model is larger (or smaller depending on the
details of the model)
than the sum of the corresponding expected log-partition functions on
$N_1$ and $N_2$-node models. This super(sub)-additivity property is used
to argue the existence of the (thermodynamic) limit of the expected
log-partition function scaled by~$N$. From this property the existence
of the scaling
limits for the ground states (optimization problems described above)
can also be shown by taking a limit as positive temperature approaches
zero temperature.
%%PT: rewrote the sentence below
In \cite{FranzLeone}, the method was used to prove the scaling limit
of log-partition functions corresponding to random K-SAT model for even $K$,
and also for the so-called Viana--Bray models with random
symmetric Hamiltonian functions.
The case of odd $K$ was apparently resolved later using the same
method \cite{FranzMontanariPrivateCommunication}.\vspace*{9pt}
%%PT: suppressing the following
%Additionally, the interpolation method was also used to establish the
%validity of replica symmetry and replica %symmetry breaking bounds on
%the free %energy limits \cite{FranzLeone,PanchenkoTalagrand},
%but %these bounds will not be discussed in this paper.

\textit{Results and technical contributions.} The goal of the present
work is to simplify and extend the applicability of the interpolation
method, and we do this in several important ways.
First, we extend the interpolation method to a variety of models on
Erd{\"o}s--R\'{e}nyi graphs
not considered before. Specifically, we consider independent set,
MAX-CUT, Ising, graph coloring (henceforth referred to as coloring),
K-SAT and Not-All-Equal K-SAT (NAE-K-SAT) models. The coloring model,
in particular, is of special interest as it is the first nonbinary
model to which interpolation method is applied.

Second, we provide a simpler and a more combinatorial interpolation
scheme as well as analysis. Moreover, we treat
the zero temperature case (optimization problem) directly and
separately from the case of the log-partition function, and again the
analysis turns out to be substantially simpler. As a result, we prove
the existence of the limit of the appropriately rescaled value of the
optimization problems in these models, including the independent set
problem, thus resolving the open problem stated earlier. %in

Third, we extend the above results to the case of random regular graphs
(and hypergraph ensembles, depending on the model). The case of random
regular graphs has been considered before by Franz, Leone and
Toninelli \cite{FranzLeoneToninelliRegular} for the K-SAT and
Viana--Bray models
with an \textit{even} number of variables per clause,
and Montanari \cite{MontanariLDPCInterpolation} in the context of
bounds on the performance of \textit{low density parity check} (LDPC)
codes. In fact, both papers consider
general degree distribution models. The second of these papers
introduces a %% MB: Removed the "more complicated"
multi-phase interpolation scheme. In this paper we consider a
modification of the interpolation scheme used in
\cite{FranzLeoneToninelliRegular} and apply it to the same six models we are
focusing in the
case of Erd{\"o}s--R\'{e}nyi graph.

Finally, we prove the large deviation principle for the satisfiability
property for coloring, K-SAT and NAE-K-SAT models on Erd{\"o}s--R\'{e}nyi graph in the following
sense. A well-known satisfiability conjecture \cite{Friedgut} states
that for each of these models there exists a (model dependent) critical
value $c^*$ such that
for every $\varepsilon>0$,
when the number of edges (or clauses for a SAT-type problem) is at most
$(c^*-\varepsilon)N$, the model is colorable (satisfiable) w.h.p., and when
it is at least $(c^*+\varepsilon)N$, it is not colorable (not
satisfiable) w.h.p. as $N\rightarrow\infty$. Friedgut
\cite{Friedgut} came close to proving
this conjecture by showing that these models exhibit sharp phase
transition: there exists a sequence $c^*_N$ such that for every
$\varepsilon$, the model is colorable (satisfiable) w.h.p. as
$N\rightarrow\infty$
when the number of edges (clauses) is at most $(c_N^*-\varepsilon)N$, and
is not colorable (satisfiable) w.h.p.
when the number of edges (clauses) is at least $(c_N^*+\varepsilon)N$.
It is also reasonable to conjecture
(which in fact is known to be true in the case
$K=2$), that not only the satisfiability conjecture is valid, but,
moreover, the probability of satisfiability $p(c,N)$ decays to zero
exponentially fast
when $c>c^*$.

In this paper we show that for these three models, namely coloring,
K-SAT and NAE-K-SAT, the limit
$r(c)\triangleq\lim_{N\to\infty}N^{-1}\log p(c,N)$ exists for every
$c$. Namely, while we do not prove the satisfiability conjecture and
the exponential
rate of convergence to zero of the satisfiability probability above the
critical threshold,
we do prove that if the convergence to zero occurs exponentially fast,
it does so at a well-defined rate $r(c)$.
Assuming the validity of the satisfiability conjecture and the
exponential rate of decay to zero above $c^*$, our result
implies that $r(c)=0$ when $c<c^*$ and $r(c)<0$ when $c>c^*$. Moreover,
we show that our results would imply the satisfiability conjecture, if
one could strengthen
Friedgut's result as follows: for every $\varepsilon>0$,
$p(c_N^*+\varepsilon,N)$ converges to zero exponentially fast, where
$c_N^*$ is the same sequence as
in Friedgut's theorem.\vspace*{9pt}

\textit{Organization of the paper.} The remainder of the paper is
organized as follows. In the following section we introduce the sparse
random (Erd{\"o}s--R\'{e}nyi) and random regular (hyper)-graphs and
introduce various combinatorial models of interest. Our main results
are stated in Section \ref{sectionMainResults}.
The proofs for the case of Erd{\"o}s--R\'{e}nyi graphs are presented in
Section \ref{sectionProofs} for results related to combinatorial optimization,
and in Section \ref{sectionProofsLogPartition} for results related to
the log-partition function. The proofs of results for random regular
graphs are presented in Section \ref{sectionProofsReg}. Several
auxiliary technical results are established in the Appendices \ref{appA} and \ref{appB}.
In particular we state and prove a simple modification of a classical
super-additivity
theorem: if a sequence is nearly super-additive, it has a limit after
an appropriate normalization.\vspace*{9pt}

\textit{Notations.} We close this section with several notational
conventions. $\R(\R_+)$ denotes the set of (nonnegative) real
values, and
$\Z(\Z_+)$ denotes the set of (nonnegative) integer values.\vadjust{\goodbreak} The
$\log$ function is assumed to be with a natural base.
As before, $[N]$ denotes the set of integers $\{1,\ldots, N\}$.
$O(\cdot)$ stands for standard order of magnitude notation.
Specifically, given two positive functions $f(N),g(N)$ defined on $N\in
\Z_+$,
$f=O(g)$ means $\sup_N f(N)/g(N)<\infty$. Also $f=o(g)$ means $\lim
_{N\rightarrow\infty} f(N)/g(N)=0$.
Throughout the paper, we treat $[N]$ as a set of nodes, and we consider
splitting this into two sets of nodes, namely $[N_1]=\{1,\ldots, N_1\}$
and $\{N_1+1,\ldots, N\}$. For symmetry, with some abuse of notation,
it is convenient to denote the second set by $[N_2]$
where $N_2=N-N_1$. $\Delta$ denotes the set-theoretic symmetric difference.
$\Bi(N,\theta)$ denotes the binomial distribution with $N$ trials and
success probability $\theta$.
$\Pois(c)$ denotes a Poisson distribution with parameter $c$,
$\stackrel{d}{=}$~stands for equality in distribution.
A sequence of random variables $X_N$ is said to converge to a random
variable $X$ with high probability (w.h.p.) if for every $\varepsilon>0$,
$\lim_{N\rightarrow\infty}\pr(|X_N-X|>\varepsilon)=0$. This is the
usual convergence in probability.

%s2 #&#
\section{Sparse random hypergraphs} %(diluted spin glass models)
\label{sectionModels} Given a set of nodes $[N]$ and a positive
integer $K$, a directed hyperedge is any ordered set of nodes
$(i_1,\ldots,i_K)\in[N]^K$. An undirected hyperedge is an unordered set
of $K$ not necessarily distinct nodes $i_1,\ldots,i_K\in[N]$. A
directed (undirected) $K$-uniform hypergraph on the node set $[N]$ is a
pair $([N],E)$, where $E$ is a set of directed (undirected)
$K$-hyperedges $E=\{e_1,\ldots,e_{|E|}\}$. Here uniformity corresponds
to the fact that every hyperedge has precisely $K$ nodes. A hypergraph
is called simple if the nodes within each hyperedge $e_m, 1\le
m\le|E|$,
are distinct and all the hyperedges are distinct. A (directed or
undirected) hypergraph is called $r$-regular if each node $i\in[N]$
appears in exactly $r$ hyperedges. The necessary condition for such a
hypergraph to exist is $Nr/K\in\Z_+$. A degree $\Delta_i=\Delta_i(\G)$
of a node $i$ is the number of hyperedges containing $i$. A (partial)
matching is a set of hyperedges such that each node belongs to at most
one hyperedge. A matching is perfect if every node of the graph belongs
to exactly one hyperedge. In this paper we use the terms hypergraph and
graph (hyperedge and edge) interchangeably.

In order to address a variety of models in a unified way, we introduce
two random directed hypergraph models, namely the Erd{\"o}s--R\'{e}nyi
random graph model $\G(N,M), M\in\Z_+$, and the random regular graph
$\G(N,r), r\in\Z_+$. These two graph models, each consisting of $N$
nodes, are described as follows. The first $\G(N,M,K)$ is obtained by
selecting $M$ directed hyperedges uniformly at random with replacement
from the space of all $[N]^K$ hyperedges. A variant of this is the
\textit{simple} Erd{\"o}s--R\'{e}nyi graph also denoted for convenience
by $\G(N,M)$, which is obtained by selecting $M$ edges uniformly at
random without replacement from the set of all undirected hyperedges
each consisting of \textit{distinct} $K$ nodes. In this paper we will
consider exclusively the case when $M=\lfloor cN\rfloor$, and $c$ is a
positive constant which does not grow with\vadjust{\goodbreak} $N$. In this case the
probability distribution of the degree of a typical node is
$\Pois(cK)+O(1/N)$. For this reason we will also call it a
\textit{sparse} random Erd{\"o}s--R\'{e}nyi graph. Often a sparse
random Erd{\"o}s--R\'{e}nyi graph is defined by including each
hyperedge in $[N]^K$ into the hypergraph with probability $c/N^{K-1}$,
and not including it with the remaining probability $1-c/N^{K-1}$. The
equivalence of two models is described using the notion of contiguity
and is well described in a variety of books, for example,
\cite{AlonSpencer,JansonBook}.

The second model $\G(N,r,K)$ is defined to be an $r$-regular directed
$K$-uniform hypergraph generated uniformly at random from
the space of all such graphs. We assume $Nr/K\in\Z_+$, so that the
set of such graphs is nonempty.
A simple (directed or undirected) version of $\G(N,r,K)$ is defined similarly.
In this paper we consider exclusively the case when $r$ is a constant
(as a function of $N$),
and we call $\G(N,r,K)$ a \textit{sparse} random regular graph.

%re1 #&#
\begin{remark}
The reason for considering the more general case of
hypergraphs is to capture
combinatorial models with hyperedges. For example, in the case of K-SAT
each clause contains $K\ge2$ distinct nodes that can be considered as
a hyperedge on $K$ nodes (more detail is provided below).
\end{remark}

%re2 #&#
\begin{remark}
In all models studied in this paper, except for K-SAT
and NAE-K-SAT, $K$ satisfies $K=2$. Therefore, to simplify the notation
we drop the reference to $K$ and throughout the paper use the shorter
notation $\G(N,M)$ and $\G(N,r)$ for the two random graph models.
\end{remark}

\textit{From nonsimple to simple graphs.} While it is common to work
with simple hypergraphs, for our purpose it is more convenient to
establish results for directed nonsimple hypergraphs first. It is well
known, however, that both $\G(N,M)$ and $\G(N,r)$ graphs are simple
with probability which remains at least a constant as $N\rightarrow
\infty$, as long as $c,r,K$ are constants. Since we prove statements
which hold w.h.p., our results have immediate ramification for simple
Erd{\"o}s--R\'{e}nyi and regular graphs.

It will be useful to recall the so-called configuration method of
constructing the random regular graph \cite{BollobasBook,BollobasRandReg,gallagerLDPC}. Each node $i$ is
associated with $r$ nodes denoted by $j^i_1,\ldots,j^i_r$. We obtain a
new set of $Nr$ nodes. Consider the $K$-uniform perfect matching
$e_1,\ldots,e_{Nr/K}$ generated uniformly at random on this set of
nodes. From this set of edges we generate a graph on the original $N$
nodes by projecting each edge to its representative. Namely an edge
$(i_1,\ldots,i_K)$ is created if and only if there is an edge of the
form $(j^{i_1}_{k_1},\ldots,j^{i_K}_{k_K})$ for some
$k_1,\ldots,k_K\in[r]$. The resulting graph is a random $r$-regular
(not necessarily simple) graph, which we again denote by $\G(N,r)$.
From now on when we talk about configuration graph, we have in mind the
graph just described on $Nr$ nodes. It is known \cite{JansonBook} that
with probability bounded away from zero as $N\rightarrow\infty$ the
resulting graph is in fact simple.

Given a hypergraph $\G=([N],E)$, we will consider a variety of
combinatorial structures on $\G$, which can be defined
in a unified way\vadjust{\goodbreak} using the notion of a Markov random field (MRF). The
MRF is a hypergraph $\G$ together
with an alphabet %% Mohsen: Changing definition of [q]
$\chi=\{0,1,\ldots,q-1\}$, denoted by $[q^-]$, and a set of node and
edge potentials $H_i, i\in[N], H_{e}, e\in E$.
A node potential is a function $H_i\dvtx [q^-]\rightarrow\R$, and an edge
potential is a function $H_e\dvtx [q^-]^K\rightarrow\lbrace-\infty\rbrace
\cup\R$.
Given a MRF $(\G,[q^-],H_i,H_e, i\in[N],e\in E)$ and any $x\in
[q^-]^N$, let
\[
H(x)=\sum_{i\in[N]} H_i(x_i)+
\sum_{e\in E}H_{e}(x_e),\qquad
H(\G)=\sup_{x\in[q^-]^N}H(x),
\]
where $x_e=(x_i, i\in e)$.
Namely, $H(x)$ is the value associated with a chosen assignment $x$,
and $H$ is the optimal value, or the groundstate
in the statistical physics terminology.
In many cases the node and edge potentials will be random functions
generated i.i.d.; see examples below.

Associated with a MRF is the Gibbs probability measure $\mu_{\G}$ on
the set of node values $[q^-]^N$ defined as follows.
Fix a parameter $\lambda>0$, and assign the probability mass $\mu_{\G
}(x)=\lambda^{H(x)}/Z_{\G}$
to every assignment $x\in[q^-]^N$, where $Z_{\G}=\sum_x\lambda
^{H(x)}$ is the normalizing \textit{partition function}.
Observe that $\lim_{\lambda\to\infty}(\log\lambda)^{-1} \log
Z_{\G}=H(\G)$. Sometimes one considers $\lambda=\exp(1/T)$\break
where $T$ is \textit{temperature}. The case $T=0$, namely $\lambda
=\infty$ then corresponds to the zero temperature regime,
or equivalently the optimization (groundstate) problem. We distinguish
this with a positive temperature case, namely $\lambda<\infty$.

We will consider in this paper a variety of MRF defined on sparse
random graphs $\G(N,\lfloor cN\rfloor)$ and $\G(N,r)$.
(In the statistical physics literature $x_i$ are called spin values,
and the corresponding MRF is called a
\textit{diluted spin glass model}.)
We now describe some examples of concrete and well-known MRF and show
that they fit the framework described above.\vspace*{9pt}

\textit{Independent set}. $K=2$ and $q=2$. Define $H_i(1)=1, H_i(0)=0$
for all $i\in[N]$.
Define $H_e(1,1)=-\infty, H_e(1,0)=H_e(0,1)=H_e(0,0)=0$ for every edge
$e=(i_1,i_2)$.
Then for every vector $x\in\{0,1\}^N$ we have
$H(x)=-\infty$ if there exists an edge $e_j=(i_1,i_2)$ such that
$x_{i_1}=x_{i_2}=1$ and
$H(x)=|\{i\dvtx x_i=1\}|$, otherwise. Equivalently, $H(x)$ takes finite
value only on $x$ corresponding to independent sets, and
in this case it is the cardinality of the independent set. $H(\G)$ is
the cardinality of a largest independent set. Note
that one can have many independent sets with cardinality $H(\G)$.\vspace*{9pt}

\textit{MAX-CUT}. $K=2$ and $q=2$. Define $H_i(0)=H_i(1)=0$.
Define $H_e(1,1)=H_e(0,0)=0, H_e(1,0)=H_e(0,1)=1$.
Every vector $x\in\{0,1\}^N$ partitions nodes into two subsets of
nodes taking values $0$ and $1$, respectively.
$H(x)$ is the number of edges between the two subsets. $H(\G)$ is the
largest such number, also called maximum cut size.
A more general case of this model is $q$-coloring; see below.\vspace*{9pt}

\textit{Anti-ferromagnetic Ising model}.
$K=2$ and $q=2$. Fix $\beta>0,B\in\R$. Define $H_i(0)=-B,H_i(1)=B$.
Define $H_e(1,1)=H_e(0,0)=-\beta, H_e(1,0)=\break H_e(0,1)=\beta$.
It is more common to use alphabet $\{-1,1\}$ instead of $\{0,1\}$ for
this model.
We use the latter for consistency with the remaining models. The
parameter~$B$, when it is nonzero
represents the presence of an external magnetic field.\vspace*{9pt}

\textit{$q$-coloring} $K=2$ and $q$ is arbitrary. $H_i(x)=0, \forall
x\in[q^-]$ and $H_e(x,y)=0$ if $x=y$ and $H_e(x,y)=1$ otherwise.
Therefore for every $x\in[q^-]^N, H(x)$ is the number of properly
colored edges, and $H(\G)$ is the maximum number of properly colored edges.\vspace*{9pt}

\textit{Random K-SAT}. $K\ge2$ is arbitrary, $q=2$.
$H_i=0$ for all $i\in[N]$.
The edge potentials $H_e$ are defined as follows. For each edge $e\in
E$ generate $a_e=(a_1,\ldots,a_K)$ uniformly at random from $\{0,1\}^K$,
independently for all edges.
For each edge $e$ set $H_e(a_1,\ldots,a_K)=0$ and $H_e(x)=1$ for all
other $x=(x_1,\ldots,x_K)$. Then for every $x\in\{0,1\}^N, H(x)$ is
the number
of satisfied clauses (hyperedges), and $H(\G)$ is the largest number
of satisfiable clauses. Often this model is called (random) MAX-K-SAT model.
We drop the MAX prefix in the notation.\vspace*{9pt}

\textit{NAE-K-SAT} (\textit{Not-All-Equal-K-SAT}). The setting is as above
except now we set $H_e(a_1,\ldots,a_K)=H_e(1-a_1,\ldots,1-a_K)=0$ and
$H_e(x)=1$ for all other $x$ for each $e$.

It is for the K-SAT and NAE-K-SAT models that considering directed, as
opposed to undirected, hypergraphs is convenient,
as for these models the order of nodes in edges matters.
For the remaining models, however, this is not the case.

In several examples considered above we have had only two possible
values for the edge potential $H_e$ and one value for the
node potential. Specifically, for the cases of coloring, K-SAT and
NAE-K-SAT problems, $H_e$ took only values $0$ and~$1$.
It makes sense to call instances of such problems ``satisfiable'' if
$H(\G)=|E|$; namely every edge potential takes value $1$.
In the combinatorial optimization terminology this corresponds to
finding a
proper coloring, a satisfying assignment and a NAE satisfying
assignment, respectively. We let
$p(N,M)=\pr(H(\G(N,M))=M)$ denote the probability of satisfiability
when the underlying graph is the Erd{\"o}s--R\'{e}nyi graph $\G(N,M)$.
We also let $p(N,r)=\pr(H(\G(N,r))=rNK^{-1})$ denote the
satisfiability probability for a random regular graph $\G(N,r)$.

%s3 #&#
\section{Main results}\label{sectionMainResults}
We now state our main results. Our first set of results concerns the
Erd{\"o}s--R\'{e}nyi graph $\G(N,\lfloor cN\rfloor)$.
%
%th1 #&#
\begin{theorem}\label{theoremMainResultER}
For every $c>0$, and for every one of the six models described in
Section \ref{sectionModels}, there exists (model dependent) $H(c)$
such that
%
%e1 #&#
\begin{equation}
\label{eqlimitexp} \lim_{N\rightarrow\infty} N^{-1}H \bigl(\G
\bigl(N,\lfloor cN\rfloor\bigr) \bigr)=H(c),
\end{equation}
w.h.p. Moreover, $H(c)$ is a Lipschitz continuous function with
Lipschitz constant~$1$.
It is a nondecreasing function of $c$ for MAX-CUT, coloring, K-SAT and
NAE-K-SAT models, and is a nonincreasing function of $c$
for the independent set model.

Also for every $c>0$ there exists $p(c)$ such that
%
%e2 #&#
\begin{equation}
\label{eqlimitprob} \lim_{N\rightarrow\infty} N^{-1}\log p\bigl(N,
\lfloor cN\rfloor\bigr)=p(c)
\end{equation}
for coloring, K-SAT and NAE-K-SAT models.
\end{theorem}

As a corollary, one obtains the following variant of the satisfiability
conjecture.
%
%co1 #&#
\begin{coro}\label{coroSATthreshold}
For coloring, K-SAT and NAE-K-SAT models, there exists a critical value
$c^*_H$ such that $H(c)=c$ when $c<c^*_H$ and $H(c)<c$ when $c>c^*_H$.
Similarly, there exists $c^*_p$, such that $p(c)=0$ when $c<c^*_p$ and
$p(c)<0$ when $c>c^*_p$.
\end{coro}
Namely, there exists a threshold value $c^*$ such that if $c<c^*$
there exists w.h.p. as $N\rightarrow\infty$ a \textit{nearly}
satisfiable assignment [assignment satisfying all but $o(N)$ clauses],
and if
$c>c^*$, then w.h.p. as $N\rightarrow\infty$, every assignment
violates linearly in $N$ many clauses. The interpretation for coloring
is similar.
The result above was established earlier by the second author for
randomly generated linear programming problems, using
the local weak convergence and martingale techniques
\cite{gamarnikLSAT}. It would be interesting to see if the same result
is obtainable using the interpolation method.

Can one use Corollary \ref{coroSATthreshold} to prove the
satisfiability conjecture in the precise sense?
The answer would be affirmative, provided that a stronger version of
Friedgut's result \cite{Friedgut} on the sharp thresholds for
satisfiability properties holds.
%
%co1 #&#
\begin{conj}\label{conjexpdecay}
For the coloring, K-SAT and NAE-K-SAT models there exists a sequence $M^*_N$
such that for every $\varepsilon>0$ there exists $\gamma=\gamma
(\varepsilon)$ such that
$\lim_{N\rightarrow\infty} p(N,\lfloor(1-\varepsilon)M^*_N\rfloor)=1$
and $p(N,\lfloor(1+\varepsilon)M^*_N\rfloor)=O(\exp(-\gamma N))$,
for all $N$.
\end{conj}
In contrast, Friedgut's sharp phase transition result \cite{Friedgut}
replaces the second part of this conjecture
with (a weaker) statement $\lim_{N\rightarrow\infty} p(N,\lfloor
(1+\varepsilon)M^*_N\rfloor)=0$.
Thus, we conjecture that beyond the phase transition region $M^*_N$,
not only is the model \textit{not satisfiable} w.h.p.,
but in fact the probability of satisfiability converges to zero
exponentially fast. The import of this (admittedly bold) statement is
as follows:

\textit{Conjecture} \ref{conjexpdecay} \textit{together with Theorem}
\ref{theoremMainResultER}
\textit{implies the satisfiability conjecture.}
Indeed, it suffices to show that $c_h^*$ is the satisfiability
threshold. We already know that
for every $\varepsilon>0$, $p(N,\lfloor(1+\varepsilon)c^*_hN)\rightarrow
0$, since $H((1+\varepsilon)c_h^*)<(1+\varepsilon)c_h^*$. Now, for the
other part
it suffices to show\vadjust{\goodbreak} that $\liminf_N M_N^*/N\rightarrow c_h^*$. Suppose
not, namely there exists $\varepsilon>0$ and a sequence $N_k$ such that
$(M_{N_k}^*/N_k)+\varepsilon<c_h^*$ for all~$k$. Then
$(M_{N_k}^*/N_k)+\varepsilon/2<c_h^*-\varepsilon/2$, implying that
%
%e3 #&#
\begin{equation}
\label{eqnearlySAT} {H(\G(N_k,\lfloor M^*_{N_k}+(\varepsilon/2)N_k\rfloor
))\over
M^*_{N_k}+(\varepsilon/2)N_k}\rightarrow1,
\end{equation}
w.h.p. by Corollary \ref{coroSATthreshold}.
On the other hand, since $M_N^*$ grows at most linearly with~$N$, we
may say
$M^*_{N_k}+(\varepsilon/2)N_k\ge(1+\varepsilon')M^*_{N_k}$, for some
$\varepsilon'>0$ for all $k$. By Conjecture \ref{conjexpdecay}, this implies
that $p(N_k,\lfloor M^*_{N_k}+(\varepsilon/2)N_k\rfloor)\rightarrow0$
exponentially fast in $N_k$. This in turn means that there exists a sufficiently
small $\delta>0$ such that the deletion of \textit{every} $\delta N_k$
edges (clauses) keeps the instance unsatisfiable w.h.p. Namely,
$H(\G(N_k,\lfloor M^*_{N_k}+(\varepsilon/2)N_k\rfloor))\le
M^*_{N_k}+(\varepsilon/2)N_k-\delta N_k$, w.h.p. as $k\rightarrow\infty$,
which contradicts (\ref{eqnearlySAT}).

Let us now state our results for the existence of the scaling limit for
the log-partition functions.
%
%th2 #&#
\begin{theorem}\label{theoremlogZ-MainResultER}
For every $c>0, 1\le\lambda< \infty$, and for every one of the
models described in Section \ref{sectionModels}, there exists (model
dependent) $z(c)$ such that
%
%e4 #&#
\begin{equation}
\label{eqlimitexpLogZ} \lim_{N\rightarrow\infty} N^{-1}\log Z\bigl(
\G\bigl(N,\lfloor cN\rfloor\bigr)\bigr)=z(c),
\end{equation}
w.h.p., where $z(c)$ is a Lipschitz continuous function of $c$.
Moreover, $z(c)$ is nondecreasing
for MAX-CUT, coloring, K-SAT and NAE-K-SAT models, and is a
nonincreasing function of $c$
for the independent set model.
\end{theorem}
%
%re3 #&#
\begin{remark}
The case $\lambda=1$ is actually not interesting as
it corresponds to no interactions between the nodes leading
to $Z(\G)=\prod_{i\in[N]}\lambda^{\sum_{x\in[q^-]}H_i(x)}$. In this
case the limit of $N^{-1}\log Z(\G(N,\lfloor cN\rfloor))$
exists trivially when node potentials $H_i$ are i.i.d.
For independent set, our proof holds for $\lambda<1$ as well. But,
unfortunately our proof
does not seem to extend to the case $\lambda<1$ in the other models.
For the Ising model this corresponds to the ferromagnetic case and the existence
of the limit was established in \cite{DemboMontanariIsing} using a
local analysis technique. The usage of local techniques is also
discussed in \cite{DemboMontanariSurvey} and \cite{DemboMontanariSun}.
%Need to assume that $\lambda< \infty$.
Finally, we remark that the proof assumes the finiteness of $\lambda$.
In fact, Coja-Oghlan observed \cite{AminPC} that if the above theorem
could suitably be extended (addressing the case of when the number of
solutions might be zero), to include the case of $\lambda=\infty$,
then the satisfiability conjecture would follow.
\end{remark}

We now turn to our results on random regular graphs.
%
%th3 #&#
\begin{theorem}\label{theoremMainResultRegular}
For every $r\in\Z_+$, and for all of the models described in the
previous section, there exists (model dependent) $H(r)$ such that
\[
\lim_{N\to\infty, N\in r^{-1}K\Z_+} N^{-1}H\bigl(\G(N,r)\bigr)=H(r)
\qquad\mbox{w.h.p.}
\]
\end{theorem}
Note, that in the statement of the theorem we take limits along
subsequence $N$ such that $NrK^{-1}$ is an integer, so that the resulting
random hypergraph is well-defined. Unlike the case of Erd{\"o}s--R\'{e}nyi graph, we were unable to prove the existence of the large deviation
rate
\[
\lim_{N\rightarrow\infty, N\in r^{-1}K\Z_+} N^{-1}\log p(N,r)
\]
for the coloring, K-SAT and NAE-K-SAT problems and leave those as open
questions.

Finally, we state our results for the log-partition function limits for
random regular graphs.
%
%th4 #&#
\begin{theorem}\label{theoremlogZ-MainResultRegular}
For every $r\in\Z_+, 1 \le\lambda< \infty$, and for every one of
the six models described in the previous section, there exists (model
dependent) $z(r)$ such that w.h.p., we have
%
%e5 #&#
\begin{equation}
\label{eqlimitexpLogZreg} \lim_{N\rightarrow\infty} N^{-1}\log Z
\bigl(\G(N,r)\bigr)=z(r).
\end{equation}
\end{theorem}

%s4 #&#
\section{Proofs: Optimization problems in Erd{\"o}s--R\'{e}nyi
graphs}\label{sectionProofs}
The following simple observation will be useful throughout the paper.
Given two hypergraphs $\G_i=([N],E_i), i=1,2$ on the same set of nodes $[N]$
for each one of the six models in Section \ref{sectionModels},
%
%e6 #&#
\begin{equation}
\label{eqobservation} \bigl|H(\G_1)-H(\G_2)\bigr|=L|E_1
\Delta E_2|,
\end{equation}
where we can take $L=1$ for all the models except Ising, and we can
take $L=\beta$ for the Ising model.
This follows from the fact that adding (deleting) an edge to (from) a
graph changes the value of $H$ by at most $1$ for all models except for
the Ising model, where the constant is $\beta$.

Our main technical result leading to the proof of Theorem
\ref{theoremMainResultER} is as follows.

%th5 #&#
\begin{theorem}\label{theoremsubadditivity}
For every $1\le N_1, N_2\le N-1$ such that $N_1+N_2=N$, and all models
%
%e7 #&#
\begin{equation}\label{eqsubadditive}
\E\bigl[H\bigl(\G\bigl(N,\lfloor cN\rfloor\bigr)\bigr)\bigr]\ge\E\bigl
[H\bigl(
\G(N_1,\rM_1)\bigr)\bigr]+\E\bigl[H\bigl(\G
(N_2,\rM_2)\bigr)\bigr],
\end{equation}
where $\rM_1\distr\Bi(\lfloor cN\rfloor, N_1/N)$ and $\rM
_2\triangleq\lfloor cN\rfloor-\rM_1\distr\Bi(\lfloor cN\rfloor, N_2/N)$.

Additionally, for the same choice of $\rM_j$ as
above and for coloring, K-SAT and NAE-K-SAT models,
%
%e8 #&#
\begin{equation}\label{eqsubadditiveProb}
p\bigl(N,\lfloor cN\rfloor\bigr) \ge\pr\bigl(H \bigl(\G(N_1,
\rM_1)\oplus\G(N_2,\rM_2) \bigr)=
\rM_1+\rM_2 \bigr),
\end{equation}
and $\G_1\oplus\G_2$ denotes a disjoint union of graphs $\G_1,\G_2$.
\end{theorem}

%re4 #&#
\begin{remark}
The event $H (\G(N_1,\rM_1)\oplus\G(N_2,\rM _2) )=\rM_1+\rM_2$
considered above corresponds to the event that both random graphs are
satisfiable (colorable) instances. The randomness of choices of edges
within each graph is assumed to be independent, but the number of edges
$\rM_j$ are dependent since they sum to $\lfloor cN\rfloor$. Because of
this coupling, it is not the case that
\begin{eqnarray*}
&&\pr\bigl(H \bigl(\G(N_1,\rM_1)\oplus\G(N_2,
\rM_2) \bigr)=\rM_1+\rM_2 \bigr)
\\
&&\qquad=\pr\bigl(H\bigl(\G(N_1,\rM_1)\bigr)=\rM_1
\bigr)\pr\bigl(H\bigl(\G(N_2,\rM_2)\bigr)=
\rM_2 \bigr).
\end{eqnarray*}
\end{remark}

Let us first show that Theorem \ref{theoremsubadditivity} implies
Theorem \ref{theoremMainResultER}.
\begin{pf*}{Proof of Theorem \ref{theoremMainResultER}}
Since $\rM_j$ have binomial distribution, we have $\E[|\rM_j-\lfloor
cN_j\rfloor|]=O(\sqrt{N})$. This together with observation
(\ref{eqobservation}) and Theorem~\ref{theoremsubadditivity} implies
\[
\E\bigl[H\bigl(\G\bigl(N,\lfloor cN\rfloor\bigr)\bigr)\bigr]\ge\E\bigl
[H\bigl(\G
\bigl(N_1,\lfloor cN_1\rfloor\bigr)\bigr)\bigr] +\E\bigl[H
\bigl(\G\bigl(N_2,\lfloor cN_2\rfloor\bigr)\bigr)\bigr]-O(
\sqrt{N}).
\]
Namely the sequence $\E[H(\G(N,\lfloor cN\rfloor))]$ is ``nearly''
super-additive, short of the $O(\sqrt{N})$ correction term.
Now we use Proposition \ref{propsubadditivityAppendix} in Appendix \ref{appB}
for the case $\alpha=1/2$ to conclude that the limit
$
\lim_{N\rightarrow\infty} N^{-1}\E[H(\G(N,\lfloor cN\rfloor
))]\triangleq H(c)
$
exists.

Showing that this also implies convergence of
$H(\G(N,\lfloor cN\rfloor))/N$ to $H(c)$
w.h.p. can be done using standard concentration results
\cite{JansonBook}, and we skip the details.
It remains to show that $H(c)$ is a nondecreasing continuous function
for MAX-CUT, coloring, K-SAT
and NAE-K-SAT problems and is nonincreasing for the independent set problem.
For the MAX-CUT, coloring, K-SAT
and NAE-K-SAT problems, we have
\[
\E\bigl[H\bigl(\G(N,M_1)\bigr)\bigr]\le\E\bigl[H\bigl(
\G(N,M_2)\bigr)\bigr],
\]
when $M_1\le M_2$; adding hyperedges can only increase the objective
value since the edge potentials are nonnegative.
For the Independent set problem on the contrary
\[
\E\bigl[H\bigl(\G(N,M_1)\bigr)\bigr]\ge\E\bigl[H\bigl(
\G(N,M_2)\bigr)\bigr]
\]
holds.
The Lipschitz continuity follows from (\ref{eqobservation}) which implies
\[
\bigl|\E\bigl[H\bigl(\G(N,M_1)\bigr)\bigr]-\E\bigl[H\bigl(
\G(N,M_2)\bigr)\bigr]\bigr|=L|M_1-M_2|
\]
with $L=\beta$ for the Ising model, and $L=1$ for the remaining models.
This concludes the proof of (\ref{eqlimitexp}).

We now turn to the proof of (\ref{eqlimitprob}) and use (\ref
{eqsubadditiveProb}) for this goal.
Our main goal is establishing the following superadditivity property:
%
%pr1 #&#
\begin{prop}\label{propsuperadProb}
There exist $0<\alpha<1$ such that for all $N_1,N_2$ such that $N=N_1+N_2$
%
%e9 #&#
\begin{equation}\label{eqsuperadProb}
\log p\bigl(N,\lfloor cN\rfloor\bigr)\ge\log p\bigl(N_1,\lfloor
cN_1\rfloor\bigr)+\log p\bigl(N_2,\lfloor
cN_2\rfloor\bigr)-O\bigl(N^\alpha\bigr).
\end{equation}
\end{prop}
\begin{pf*}{Proof of Proposition \ref{propsuperadProb}}
Fix any $1/2<\nu<1$. First we assume $N_1\le N^{\nu}$.
Let $\rM_j$ be as in Theorem \ref{theoremsubadditivity}. We have
\begin{eqnarray*}
&&
\pr\bigl(H \bigl(\G(N_1,\rM_1)\oplus\G(N_2,
\rM_2) \bigr)=\rM_1+\rM_2 \bigr)
\\
&&\qquad
\ge p\bigl(N_2,\lfloor cN_2\rfloor\bigr) p
\bigl(N_1,\lfloor cN\rfloor-\lfloor cN_2\rfloor\bigr)\pr
\bigl(\rM_2=\lfloor cN_2\rfloor\bigr).
\end{eqnarray*}
We have
\[
\pr\bigl(\rM_2=\lfloor cN_2\rfloor\bigr)=\pmatrix{\lfloor cN
\rfloor\cr\lfloor cN_2\rfloor} (N_2/N)^{\lfloor cN_2\rfloor
}(N_1/N)^{\lfloor cN\rfloor-\lfloor
cN_2\rfloor}.
\]
From our assumption $N_1\le N^\nu$ it follows that
\begin{eqnarray*}
(N_2/N)^{\lfloor cN_2\rfloor}&\ge&\bigl(1-N^{\nu-1}
\bigr)^{O(N)}=\exp\bigl(-O\bigl(N^\nu\bigr)\bigr),
\\
(N_1/N)^{\lfloor cN\rfloor-\lfloor cN_2\rfloor}&\ge&(1/N)^{O(N_1)}\ge
\exp\bigl(-O
\bigl(N^\nu\log N\bigr)\bigr).
\end{eqnarray*}
It then follows
\[
\pr\bigl(\rM_2=\lfloor cN_2\rfloor\bigr)\ge\exp\bigl(-O
\bigl(N^\nu\log N\bigr)\bigr).
\]
Now we claim the following crude bound for every deterministic
$m$ and every one of the three models under the consideration.
\[
p(N,m+1)\ge O(1/N)p(N,m).
\]
Indeed, for the K-SAT model, conditional on the event that $H(\G
(N,m))=m$, the probability that $H(\G(N,m+1))=m+1$ is at least
$1-1/2^K$.
We obtain thus a bound which is even stronger than claimed,
\[
p(N,m+1)\ge\bigl(1-1/2^K\bigr)p(N,m)=O\bigl(p(N,m)\bigr).
\]
The proof for the NAE-K-SAT is similar.
For the coloring problem observe that this conditional probability is
at least $(1-1/N)(2(N-1)/N^2)=O(1/N)$ since with probability $1-1/N$
the new edge chooses
different nodes, and with probability at least $2(N-1)/N^2$ the new
edge does not violate a given coloring (with equality achieved
only when $q=2$, and two coloring classes having cardinalities $1$ and
$N-1$). The claim follows.

Now since $\lfloor cN\rfloor-\lfloor cN_2\rfloor\le\lfloor
cN_1\rfloor+1$, the claim implies
\[
p\bigl(N_1,\lfloor cN\rfloor-\lfloor cN_2\rfloor\bigr)
\ge O(1/N)p\bigl(N_1,\lfloor cN_1\rfloor\bigr).
\]
Combining our estimates we obtain
\begin{eqnarray*}
&&
\pr\bigl(H \bigl(\G(N_1,\rM_1)\oplus\G(N_2,
\rM_2) \bigr)=\rM_1+\rM_2 \bigr)
\\
&&\qquad
\ge p\bigl(N_1,\lfloor cN_1\rfloor\bigr)p
\bigl(N_2,\lfloor cN_2\rfloor\bigr)O(1/N)\exp\bigl(-O
\bigl(N^\nu\log N\bigr)\bigr).
\end{eqnarray*}
After taking logarithm of both sides we obtain (\ref{eqsuperadProb})
from (\ref{eqsubadditiveProb}).\eject

The case $N_2\le N^{\nu}$ is considered similarly.
We now turn to a more difficult case $N_j>N^{\nu}, j=1,2$.

First we state the following lemma (proved in Appendix \ref{appA}) for the three
models of interest (coloring, K-SAT, NAE-K-SAT).
%
%le1 #&#
\begin{lemma}\label{lemmapN,m+1>CpN,m} The following holds for
coloring, K-SAT, NAE-K-SAT models for all $N,M,m$ and $0<\delta<1/2$:
%
%e10 #&#
\begin{equation}\label{eqpN,m+1>CpN,m}
p(N,M+m)\ge\delta^m p(N,M)-(2\delta)^{M+1}\exp\bigl(H(
\delta)N+o(N)\bigr),
\end{equation}
where $H(\delta)=-\delta\log\delta-(1-\delta)\log(1-\delta)$ is
the entropy function.
\end{lemma}

We now prove (\ref{eqlimitprob}). Fix $h\in(1/2,\nu)$.
We have from (\ref{eqsubadditiveProb}),
\begin{eqnarray*}
&&
p\bigl(N,\lfloor cN\rfloor\bigr)\\
&&\qquad\ge\pr\bigl(H \bigl(\G(N_1,
\rM_1)\oplus\G(N_2,\rM_2) \bigr)=
\rM_1+\rM_2 \bigr)
\\
&&\qquad\ge\sum_{cN_1-N^h\le m_1\le cN_1+N^h, m_2=\lfloor cN\rfloor-m_1}
p(N_1,m_1)p(N_2,m_2)
\pr(\rM_1=m_1).
\end{eqnarray*}
Note that $cN_1-N^h\le m_1\le cN_1+N^h$ implies $cN_2-N^h-1\le m_2\le cN_2+N^h$.
Applying Lemma \ref{lemmapN,m+1>CpN,m} we further obtain for the
relevant range of $m_j$ that
\begin{eqnarray*}
&&
p(N_j,m_j)\\
&&\qquad\ge\delta^{(m_j-\lfloor cN_j\rfloor)^+ } p
\bigl(N_j,\lfloor cN_j\rfloor\bigr)-(2
\delta)^{cN_j}\exp\bigl(H(\delta)N_j+o(N_j)\bigr)
\\
&&\qquad\ge\delta^{N^h+1} p\bigl(N_j,\lfloor cN_j\rfloor
\bigr)-(2\delta)^{cN_j}\exp\bigl(H(\delta)N_j+o(N_j)
\bigr)
\\
&&\qquad\ge\delta^{N^h+1}p\bigl(N_j,\lfloor cN_j\rfloor
\bigr) \biggl(1-2^{cN_j}\delta^{cN_j-N^h-1}\biggl(1-\frac
{1}{q}
\biggr)^{-cN_j}e^{H(\delta)N_j+o(N_j)} \biggr),
\end{eqnarray*}
where we have used a simple bound $p(N_j,\lfloor cN_j\rfloor)\ge
(1-1/q)^{cN_j}$. Now let us take $\delta$
so that
%
%e11 #&#
\begin{equation}
\label{eqbetadelta} \beta(\delta)\triangleq\bigl(2\delta(1-1/q)
\bigr)^{-c}\exp\bigl(H(\delta)\bigr)<1.
\end{equation}
Then using the assumptions $N_j\ge N^{\nu}$ and $h<\nu$ we obtain
\[
p(N_j,m_j)\ge\delta^{N^h+1}p\bigl(N_j,
\lfloor cN_j\rfloor\bigr) \bigl(1-\bigl(\beta(\delta)
\bigr)^{O(N^\nu)} \bigr).
\]

Combining we obtain
\begin{eqnarray*}
p\bigl(N,\lfloor cN\rfloor\bigr)&\ge&\pr\bigl(cN_1-N^h\le
\rM_1\le cN_1+N^h\bigr)
\\
&&{}\times\prod_{j=1,2}\delta^{O(N^h)} p
\bigl(N_j,\lfloor cN_j\rfloor\bigr) \bigl(1-\bigl(\beta(
\delta)\bigr)^{O(N^\nu)} \bigr).
\end{eqnarray*}
This implies
\begin{eqnarray*}
\log p\bigl(N,\lfloor cN\rfloor\bigr)&\ge&\log\pr\bigl(cN_1-N^h
\le\rM_1\le cN_1+N^h\bigr)
\\
&&{}+N^h\log\delta+\sum_{j=1,2}\log p
\bigl(N_j,\lfloor cN_j\rfloor\bigr)\\
&&{}+\log\bigl(1-\bigl(
\beta(\delta)\bigr)^{O(N^\nu)} \bigr).
\end{eqnarray*}
Since $\rM_1\distr\Bi(\lfloor cN\rfloor, N_1/N)$ and $h>1/2$, then
\[
\bigl|\log\pr\bigl(cN_1-N^h\le\rM_1\le
cN_1+N^h\bigr) \bigr|=o(1).
\]
Since $\beta(\delta)<1$, then
\[
\log\bigl(1-\bigl(\beta(\delta)\bigr)^{O(N^\nu)} \bigr)=O \bigl(\bigl
(\beta
(\delta)\bigr)^{O(N^\nu)} \bigr)=o\bigl(N^h\bigr),
\]
where the last identity is of course a very crude estimate.
Combining, we obtain
\[
\log p\bigl(N,\lfloor cN\rfloor\bigr)\ge\sum_{j=1,2}
\log p\bigl(N_j,\lfloor cN_j\rfloor\bigr)+O
\bigl(N^h\bigr).
\]
The claim of Proposition \ref{propsuperadProb} is established.
\end{pf*}

Part (\ref{eqlimitprob}) of Theorem \ref{theoremMainResultER} then
follows from this proposition and
Proposition~\ref{propsubadditivityAppendix} from Appendix \ref{appB}.
\end{pf*}

We now turn to the proof of Theorem \ref{theoremsubadditivity} and,
in particular, introduce the interpolation construction.
\begin{pf*}{Proof of Theorem \ref{theoremsubadditivity}}
We begin by constructing a sequence of graphs interpolating between $\G
(N,\lfloor cN\rfloor)$ and a disjoint union of $\G(N_1,\rM_1)$ and
$\G(N_2,\lfloor cN\rfloor-\rM_1)$. Given $N, N_1, N_2$ s.t. $N_1+N_2=N$
and any $0\le r\le\lfloor cN\rfloor$, let $\G(N,\lfloor cN\rfloor,r)$
be the random graph on nodes $[N]$ obtained as follows. It contains
precisely $\lfloor cN\rfloor$ hyperedges. The first $r$ hyperedges
$e_1,\ldots,e_r$ are selected u.a.r. from all the possible directed
hyperedges [namely they are generated as hyperedges of $\G(N,\lfloor
cN\rfloor)$]. The remaining $\lfloor cN\rfloor-r$ hyperedges
$e_{r+1},\ldots,e_{\lfloor cN\rfloor}$ are generated as follows. For
each $j=r+1,\ldots,\lfloor cN\rfloor$, with probability $N_1/N$, $e_j$
is generated independently u.a.r. from all the possible hyperedges on
nodes $[N_1]$, and with probability $N_2/N$, it is generated u.a.r.
from all the possible hyperedges on nodes $[N_2]$ (\mbox{$=$}$\{N_1+1,\ldots,
N\}$). The choice of node and edge potentials $H_v,H_e$ is done exactly
according to the corresponding model, as for the case of graphs
$\G(N,\lfloor cN\rfloor)$. Observe that when $r=\lfloor cN\rfloor$,
$\G(N,\lfloor cN\rfloor,r)=\G(N,\lfloor cN\rfloor)$, and when $r=0$,
$\G(N,\lfloor cN\rfloor,r)$ is a disjoint union of graphs
$\G(N_1,\rM_1), \G(N_2,\rM_2)$, conditioned on $\rM_1+\rM_2=\lfloor
cN\rfloor$, where $\rM_j\distr \Bi(\lfloor cN\rfloor, N_j/N)$.

%pr2 #&#
\begin{prop}\label{propsubadditivity}
For every $r=1,\ldots,\lfloor cN\rfloor$,
\[
\E\bigl[H\bigl(\G\bigl(N,\lfloor cN\rfloor,r\bigr)\bigr)\bigr] \ge\E
\bigl[H
\bigl(\G\bigl(N,\lfloor cN\rfloor,r-1\bigr)\bigr)\bigr]. %
\]
Also for coloring, K-SAT and NAE-K-SAT models,
\[
\pr\bigl(H \bigl(\G\bigl(N,\lfloor cN\rfloor,r\bigr) \bigr)=\lfloor
cN\rfloor
\bigr) \ge\pr\bigl(H \bigl(\G\bigl(N,\lfloor cN\rfloor,r-1\bigr) \bigr
)=\lfloor
cN\rfloor\bigr). %\label{eqpartialsubadditivityProb}.
\]
\end{prop}
Let us first show how
Theorem \ref{theoremsubadditivity} follows from this proposition.
Observe that for a disjoint union of two deterministic graphs $\G=\G
_1+\G_2$, with $\G=(V,E), \G_1=(V_1,E_1),\G_2=(V_2,E_2)$,
we always have $H(\G)=H(\G_1)+H(\G_2)$. Claim (\ref
{eqsubadditive}) then follows. Caim (\ref{eqsubadditiveProb})
follows immediately from the interpolation construction by comparing
the cases $r=0$ and $r=\lfloor cN\rfloor$.
\end{pf*}

\begin{pf*}{Proof of Proposition \ref{propsubadditivity}}
Observe that $\G(N,\lfloor cN\rfloor,r-1)$ is obtained from $\G
(N,\lfloor cN\rfloor,r)$ by deleting a hyperedge chosen u.a.r.
independently from $r$ hyperedges $e_1,\ldots,e_r$
and adding a hyperedge either to nodes $[N_1]$ or to $[N_2]$ with
probabilities $N_1/N$ and $N_2/N$, respectively. Let $\G_0$ be the
graph obtained
after deleting but before adding a hyperedge. For the case of K-SAT and
NAE-K-SAT (two models with random edge potentials), assume that $\G_0$
also encodes the underlying edge potentials of the instance.
For the case of coloring, K-SAT, NAE-K-SAT, note that the maximum value
that $H$ can achieve for the graph $\G_0$ is $\lfloor cN\rfloor-1$
since exactly one hyperedge was deleted.
We will establish a stronger result: conditional on \textit{any}
realization of the graph $\G_0$ (and random potentials), we
claim that
%
%e12 #&#
\begin{equation}\label{eqconditionalsubadditivity}
\E\bigl[H\bigl(\G\bigl(N,\lfloor cN\rfloor,r\bigr)\bigr)|\G_0\bigr]
\ge\E\bigl[H\bigl(\G\bigl(N,\lfloor cN\rfloor,r-1\bigr)\bigr)|\G_0
\bigr]
\end{equation}
and
%
%e13 #&#
\begin{eqnarray}\label{eqconditionalsubadditivityProb}
&&
\pr\bigl(H\bigl(\G\bigl(N,\lfloor cN\rfloor,r\bigr)\bigr)=\lfloor
cN\rfloor|
\G_0\bigr)\nonumber\\[-8pt]\\[-8pt]
&&\qquad\ge\pr\bigl(H\bigl(\G\bigl(N,\lfloor cN\rfloor,r-1\bigr)
\bigr)=\lfloor cN\rfloor|\G_0\bigr)\nonumber
\end{eqnarray}
for coloring, K-SAT, NAE-K-SAT. Proposition then follows immediately
from these claims by averaging over $\G_0$.
Observe that \textit{conditional on any realization $\G_0$}, $\G
(N,\lfloor cN\rfloor,r)$ is obtained from $\G_0$
by adding a hyperedge to $[N]$ u.a.r. That is the generation of this hyperedge
is independent from the randomness of $\G_0$. Similarly, conditional
on any realization $\G_0$, $\G(N,\lfloor cN\rfloor,r-1)$ is obtained
from $\G_0$ by adding a hyperedge to
$[N_1]$ or $[N_2]$ u.a.r. with probabilities $N_1/N$ and $N_2/N$, respectively.

We now prove properties (\ref{eqconditionalsubadditivity}) and (\ref
{eqconditionalsubadditivityProb}) for each of the six models.\vspace*{9pt}

\textit{Independent sets}.
Let $O^*\subset[N]$ be the set of nodes which belong to \textit{every}
largest independent set in $\G_0$. Namely if $I\subset[N]$ is an i.s.
such that
$|I|=H(\G_0)$, then $O^*\subset I$. We note that $O^*$ can be empty.
Then for every edge $e=(i,k)$, $H(\G_0+e)=H(\G_0)-1$ if $i,k\in O^*$ and
$H(\G_0+e)=H(\G_0)$ if either $i\notin O^*$ or $k\notin O^*$. Here
$\G_0+e$ denotes a graph obtained from $\G_0$ by adding $e$.
When the edge $e$ is generated u.a.r. from the all possible edges, we
then obtain
$
\E[H(\G_0+e)|\G_0]-H(\G_0)=- ({|O^*|\over N} )^2.
$
Therefore, $\E[H(\G(N,\lfloor cN\rfloor,r))|\G_0]-H(\G_0)=-
({|O^*|\over N} )^2$.
By a similar argument
\begin{eqnarray*}
&&\E\bigl[H\bigl(\G\bigl(N,\lfloor cN\rfloor,r-1\bigr)\bigr)|\G_0\bigr]-H(
\G_0)
\\
&&\qquad=-{N_1\over N} \biggl({|O^*\cap[N_1]| \over N_1} \biggr)^2
-{N_2\over
N} \biggl({|O^*\cap[N_2]| \over N_2} \biggr)^2
\\
&&\qquad\le- \biggl({N_1\over N}{|O^*\cap[N_1]| \over N_1} +
{N_2\over
N}{|O^*\cap[N_2]| \over N_2} \biggr)^2
\\
&&\qquad=- \biggl({|O^*|\over N} \biggr)^2=\E\bigl[H(\G\bigl(N,
\lfloor cN\rfloor,r\bigr)|\G_0\bigr]-H(\G_0),
\end{eqnarray*}
and (\ref{eqconditionalsubadditivity}) is established.\vspace*{9pt}

\textit{MAX-CUT}.
Given $\G_0$, let $\mathcal{C}^*\subset\{0,1\}^{[N]}$ be the set of
optimal solutions. Namely $H(x)=H(\G_0), \forall x\in\mathcal{C}^*$
and $H(x)<H(\G_0)$ otherwise.
Introduce an equivalency relationship $\sim$ on $[N]$. Given $i,k\in
[N]$, define $i\sim k$ if \textit{for every} $x\in\mathcal{C}^*, x_i=x_k$.
Namely, in every optimal cut, nodes $i$ and $k$ have the same value.
Let $O^*_j\subset[N], 1\le j\le J$, be the corresponding equivalency
classes. Given any edge $e=(i,k)$, observe that
$H(\G_0+e)=H(\G_0)$ if $i\sim k$ and $H(\G_0+e)=H(\G_0)+1$
otherwise. Thus
\[
\E\bigl[H(\G\bigl(N,\lfloor cN\rfloor,r\bigr)|\G_0\bigr]-H(
\G_0)=1-\sum_{1\le j\le
J} \biggl(
{|O_j^*|\over N} \biggr)^2
\]
and
\begin{eqnarray*}
&&\E\bigl[H(\G\bigl(N,\lfloor cN\rfloor,r-1\bigr)|\G_0\bigr]-H(
\G_0)
\\
&&\qquad=1-{N_1\over N}\sum_{1\le j\le J} \biggl(
{|O_j^*\cap[N_1]| \over
N_1} \biggr)^2 -{N_2\over N}\sum
_{1\le j\le J} \biggl({|O_j^*\cap[N_2]| \over
N_2}
\biggr)^2.
\end{eqnarray*}
Using ${N_1\over N} ({|O_j^*\cap[N_1]| \over N_1} )^2+
{N_2\over N} ({|O_j^*\cap[N_2]| \over N_2} )^2 \ge
({|O_j^*|\over N} )^2$
we obtain (\ref{eqconditionalsubadditivity}).\vspace*{9pt}

\textit{Ising}.
The proof is similar to the MAX-CUT problem but is more involved due to
the presence of the magnetic field $B$.
The presence of the field means that we can no longer say that $H(\G
_0+e)=H(\G_0)+\beta$ or \mbox{$=$}$H(\G_0)$. This issue
is addressed by looking at suboptimal solutions and the implied
sequence of equivalence classes.
Thus let us define a sequence $H_0>H_1>H_2\cdots>H_M$ and an integer
$M\ge0$ as follows.
Define $H_0=H(\G_0)$. Assuming $H_{m-1}$ is defined, $m\ge1$, let
$H_m=\max H(x)$ over all solutions $x\in\{0,1\}^N$ such that $H(x)<H_{m-1}$.
Namely, $H_m$ is the next best solution after $H_{m-1}$. Define $M$ to
be the largest $m$ such that $H_m\ge H(\G_0)-2\beta$.
If this is not the case for all $m$, then we define $M\le2^N$ to be
the total number of possible values $H(x)$ (although typically the
value of $M$ will be much
smaller). Let $\mathcal{C}_m=\{x\dvtx H(x)=H_m\}, 0\le m\le M$, be the set
of solutions achieving value $H_m$. Observe that
$\mathcal{C}_m$ are disjoint sets.
For every $m\le M$ define an equivalency
relationship as follows $i\stackrel{m}\sim k$ for $i,k\in[N]$ if and
only if $x_i=x_k$ for all $x\in\mathcal{C}_0\cup\cdots\cup\mathcal{C}_m$.
Namely, nodes $i$ and
$k$ are $m$-equivalent if they take equal values in every solution
achieving value at least $H_m$. Let $O^*_{j,m}$ be the corresponding
equivalency classes
for $1\le j\le J_m$. Note that the partition $O^*_{j,m+1}$ of the nodes
$[N]$ is a refinement of the partition $O^*_{j,m}$.

%le2 #&#
\begin{lemma}\label{lemmaIsing}
Given an edge $e=(i,k)$, the following holds:
\[
H(\G+e)=\cases{ H(\G_0)+\beta, &\quad if $i\stackrel{0} {\not\sim} k$;
\cr
H_{m+1}+\beta, &\quad if $i\stackrel{m}\sim k$, but $i\stackrel{m+1} {
\not\sim} k$, for some $m\le M-1$;
\cr
H(\G_0)-\beta, &\quad if $i
\stackrel{m}\sim k$ for all $m\le M$.}
\]
\end{lemma}

\begin{pf}
The case $i\stackrel{0}{\not\sim} k$ is straightforward.
Suppose $i\stackrel{m}\sim k$, but $i\stackrel{m+1}{\not\sim} k$
for some $m\le M-1$.
For every $x\in\bigcup_{m'\le m}\mathcal{C}_{m'}$ we have for some
$m'\le m$, $H_{\G_0+e}(x)=H_{m'}-\beta\le H_0-\beta$.
Now since $i\stackrel{m+1}{\not\sim} k$ there exists $x\in\mathcal
{C}_{m+1}$ such that $x_i\ne x_k$, implying
$H_{\G_0+e}(x)=H_{m+1}+\beta\ge H_0-\beta$, where the inequality
follows since $m+1\le M$. Furthermore, for every $x\notin\bigcup_{m'\le
m}\mathcal{C}_{m'}$
we have $H_{\G_0+e}(x)\le H_{\G_0}(x)+\beta\le H_{m+1}+\beta$. We
conclude that $H_{m+1}+\beta$ is the optimal solution in this case.

On the other hand, if $i\stackrel{m}\sim k$ for all $m\le M$, then for
all $x\in\bigcup_{m\le M}\mathcal{C}_{m}$,
$H_{\G_0+e}(x)\le H(\G_0)-\beta$, with equality achieved for $x\in
\mathcal{C}_0$. For all\vspace*{1pt}
$x\notin\bigcup_{m\le M}\mathcal{C}_{m}$, we have $H_{\G_0+e}(x)\le
H_{M+1}+\beta<H_0-\beta$, and the assertion is established.
Note that if $M=2^N$, namely $M+1$ is not defined, then $\bigcup_{m\le
M}\mathcal{C}_{m}$ is the entire space of solutions $\{0,1\}^N$,
and the second part of the previous sentence is irrelevant.
\end{pf}

We now return to the proof of the proposition. Recall that if an edge
$e=(i,k)$ is added uniformly at random then
$\pr(i\stackrel{m}\sim k)=\sum_{1\le j\le J_m}
({|O_{j,m}^*|\over N} )^2$. A similar\vspace*{1pt} assertion holds for the
case $e$ is added uniformly at random to parts $[N_l], l=1,2$, with
probabilities $N_l/N$, respectively.
We obtain that
\begin{eqnarray*}
\pr\bigl(i\stackrel{0} {\not\sim} k\bigr)&=&1-\sum_{1\le j\le J_0}
\biggl({|O_{j,0}^*|\over N} \biggr)^2,
\\
\pr\bigl(i\stackrel{m}\sim k\mbox{, but } i\stackrel{m+1} {\not\sim}
k\bigr)&=&
\sum
_{1\le j\le J_m} \biggl({|O_{j,m}^*|\over N}
\biggr)^2 -\sum_{1\le j\le J_{m+1}} \biggl(
{|O_{j,m+1}^*|\over N} \biggr)^2,
\\
\pr\bigl(i\stackrel{m}\sim k, \forall m\le M\bigr)&=&\sum_{1\le j\le J_M}
\biggl({|O_{j,M}^*|\over N} \biggr)^2.
\end{eqnarray*}
Applying Lemma \ref{lemmaIsing} we obtain
\begin{eqnarray*}
&&
\E\bigl[H\bigl(\G\bigl(N,\lfloor cN\rfloor,r\bigr)\bigr)|\G_0
\bigr]\\
&&\qquad=(H+\beta) \biggl(1-\sum_{1\le j\le J_m} \biggl(
{|O_{j,0}^*|\over N} \biggr)^2 \biggr)
\\
&&\qquad\quad{}+\sum_{m=0}^{M-1}(H_{m+1}+\beta)
\biggl(\sum_{1\le j\le J_m} \biggl({|O_{j,m}^*|\over N}
\biggr)^2 -\sum_{1\le j\le J_{m+1}} \biggl(
{|O_{j,m+1}^*|\over N} \biggr)^2 \biggr)
\\
&&\qquad\quad{}+ (H-\beta)\sum_{1\le j\le J_M} \biggl(
{|O_{j,M}^*|\over
N} \biggr)^2
\\
&&\qquad=H+\beta+\sum_{0\le m\le M-1}(H_{m+1}-H_{m})
\sum_{1\le j\le
J_m} \biggl({|O_{j,m}^*|\over N}
\biggr)^2
\\
&&\qquad\quad{}+(H-H_M-2\beta)\sum_{1\le j\le J_M} \biggl(
{|O_{j,M}^*|\over N} \biggr)^2.
\end{eqnarray*}
By a similar argument and again using Lemma \ref{lemmaIsing} we obtain
\begin{eqnarray*}
&&\E\bigl[H\bigl(\G\bigl(N,\lfloor cN\rfloor,r-1\bigr)\bigr)|\G_0
\bigr]
\\
&&\qquad=H+\beta+\sum_{0\le m\le M-1}(H_{m+1}-H_{m})
\sum_{1\le j\le J_m} \sum_{l=1,2}
{N_l\over N} \biggl({|O_{j,m}^*\cap[N_l]|\over N_l} \biggr)^2
\\
&&\qquad\quad{}+(H-H_M-2\beta)\sum_{1\le j\le J_M}\sum
_{l=1,2}{N_l\over N} \biggl({N_l\over N}
{|O_{j,M}^*\cap[N_l]|\over N_l} \biggr)^2.
\end{eqnarray*}
Recall, however, that $H_{m+1}-H_m<0, m\le M-1$ and $H-H_M-2\beta\le
0$. Again using the convexity of the $g(x)=x^2$ function, we obtain the claim.\vspace*{9pt}

\textit{Coloring}. Let $\mathcal{C}^*\subset[q^-]^N$ be the set of
optimal colorings. Namely $H(x)=H(\G_0), \forall x\in\mathcal{C}^*$.
Introduce an equivalency relationship $\sim$ on the set of nodes as follows.
Given $i,k\in[N]$, define $i\sim k$ if and only if $x_i=x_k$ for every
$x\in\mathcal{C}^*$. Namely, in every optimal coloring assignments,
$i$ and $k$
receive the same color. Then for every edge $e$, $H(\G_0+e)=H(\G_0)$
if $i\sim k$ and $H(\G_0+e)=H(\G_0)+1$ otherwise. The remainder
of the proof of (\ref{eqconditionalsubadditivity}) is similar to the
one for MAX-CUT.

Now let us show (\ref{eqconditionalsubadditivityProb}). We fix graph
$\G_0$. Notice that if $\G_0$ is not colorable,
then both probabilities in (\ref{eqconditionalsubadditivityProb})
are zero, since adding edges cannot turn an uncolorable instance
into the colorable one.
Thus assume $\G_0$ is a colorable graph. Since it has
$\lfloor cN\rfloor-1$ edges it means $H(\G_0)=\lfloor cN\rfloor-1$.
Let
$O_j^*\subset[N], 1\le j\le J$ denote the $\sim$ equivalence classes,
defined by $i\sim k$ if and only if in every proper coloring
assignment $i$ and $k$ receive the same color. We obtain that
\[
\pr\bigl(H \bigl(\G\bigl(N,\lfloor cN\rfloor,r\bigr) \bigr)=\lfloor
cN\rfloor|
\G_0 \bigr)=1-\sum_{1\le j\le J} \biggl(
{|O_j^*|\over N} \biggr)^2.
\]
Similarly,
\begin{eqnarray*}
&&
\pr\bigl(H \bigl(\G\bigl(N,\lfloor cN\rfloor,r-1 \bigr)\bigr)=\lfloor
cN\rfloor
|\G_0 \bigr)
\\
&&\qquad=1-{N_1\over N}\sum_{1\le j\le J} \biggl(
{|O_j^*\cap[N_1]| \over
N_1} \biggr)^2 -{N_2\over N}\sum
_{1\le j\le J} \biggl({|O_j^*\cap[N_2]| \over N_2}
\biggr)^2.
\end{eqnarray*}
Relation (\ref{eqconditionalsubadditivityProb}) then again follows
from convexity.\vspace*{9pt}

\textit{K-SAT}. Let $\mathcal{C}^*\subset\{0,1\}^N$ be the set
of optimal assignments. Define a node $i$ (variable $x_i)$ to be
\textit{frozen} if either $x_i=0, \forall x\in\mathcal{C}^*$ or $x_i=1,
\forall x\in \mathcal{C}^*$. Namely, in every optimal assignment the
value of $i$ is always the same. Let $O^*$ be the set of frozen
variables. Let $e=(i_1,\ldots,i_K)\subset[N]$ be a hyperedge, and let
$H_e\dvtx \{ 0,1\}^K\rightarrow\{0,1\}$ be the corresponding edge potential.
Namely, for some $y_1,\ldots,y_K\in\{0,1\},
H_e(x_{i_1},\ldots,x_{i_k})=0$ if $x_{i_1}=y_1,\ldots,x_{i_K}=y_K$ and
$H_e=1$ otherwise. Consider adding $e$ with $H_e$ to the graph $\G_0$.
Note that if \mbox{$e\cap([N]\setminus O^*)\ne\varnothing$}, then
$H(\G_0+e)=H(\G_0)+1$, as in this case at least one variable in $e$ is
nonfrozen and can be adjusted to satisfy the clause. Otherwise,
suppose $e\subset O^*$, and let $x_{i_1}^*,\ldots,x_{i_K}^*\in\{0,1\} $
be the corresponding frozen values of $i_1,\ldots,i_K$. Then
$H(\G_0+e)=H(\G_0)$ if $x_{i_1}^*=y_1,\ldots,x_{i_K}^*=y_K$, and
$H(\G_0+e)=H(\G_0)+1$ otherwise. Moreover, for the random choice of
$H$, the first event $H(\G _0+e)=H(\G_0)$ occurs with probability
$1/2^K$. We conclude that
\[
\E\bigl[H \bigl(\G\bigl(N,\lfloor cN\rfloor,r \bigr)\bigr) |\G_0 \bigr]-H(
\G_0)=1-{1\over2^K} \biggl({|O^*|\over N}
\biggr)^K
\]
and for every satisfiable instance $\G_0$ (namely $H(\G_0)=\lfloor
cN\rfloor-1$),
\[
\pr\bigl(H \bigl(\G\bigl(N,\lfloor cN\rfloor,{r}\bigr) \bigr)=\lfloor
cN\rfloor
|{\G_0} \bigr) =1-{1\over2^K} \biggl(
{|O^*|\over N} \biggr)^K.
\]
Similarly,
\begin{eqnarray*}
&&\E\bigl[H\bigl(\G\bigl(N,\lfloor cN\rfloor,r-1\bigr)\bigr)|\G_0,
H_0\bigr]-H(\G_0)
\\
&&\qquad=1-{1\over2^K}{N_1\over N} \biggl({|O^*\cap[N_1]| \over N_1}
\biggr)^K -{1\over2^K}{N_2\over N} \biggl(
{|O^*\cap[N_2]| \over N_2} \biggr)^K
\end{eqnarray*}
and for every satisfiable instance $\G_0$,
\begin{eqnarray*}
&&
\pr\bigl(H \bigl(\G\bigl(N,\lfloor cN\rfloor,r-1\bigr) \bigr)=\lfloor
cN\rfloor
|{\G_0} \bigr)
\\
&&\qquad=1-{1\over2^K}{N_1\over N} \biggl({|O^*\cap[N_1]| \over N_1}
\biggr)^K -{1\over2^K}{N_2\over N} \biggl(
{|O^*\cap[N_2]| \over N_2} \biggr)^K.
\end{eqnarray*}
Using the convexity of the function $x^K$ on $x\in[0,\infty)$, we
obtain the result.\vspace*{9pt}

\textit{NAE-K-SAT}.
The idea of the proof is similar and is based on the combination of the
notions of frozen variables and equivalency classes.
Two nodes (variables) $i$~and $k$ are defined to be equivalent $i\sim
k$ if there do not exist two assignments $x,x'$ such that
$x_i=x'_i$, but $x_k\ne x'_k$, or vice verse, $x_i\ne x'_i$, but
$x_k=x'_k$. Namely, either both nodes are frozen, or
setting one of them determines the value for the other in every optimal
assignment. Let $O_j^*, 1\le j\le J$, be the set of
equivalence classes (the set of frozen variables is one of $O^*_j$).
Let $e=(i_1,\ldots,i_K)\subset[N]$
be a hyperedge added to $\G_0$, and let $H_e\dvtx \{0,1\}^K\rightarrow\{
0,1\}$
be the corresponding edge potential. We claim that if $i_1,\ldots,i_K$
are not all equivalent, then $H(\G+e)=H(\G)$.
Indeed, suppose without the loss of generality that $i_1\not\sim i_2$
and $x,x'$ are two optimal solutions such that
$x_{i_1}=x'_{i_1}, x_{i_2}\ne x_{i_2}'$. From the definition of
NAE-K-SAT model, it follows that at least one of the
two solutions $x$ and $x'$ satisfies $H_e$ as well, and the claim then
follows. Thus, $H(\G+e)=H(\G)$ only if
$i_1,\ldots,i_K$ all belong to the same equivalence class. Provided
that this indeed occurs, it is easy to see
that the probability that $H(\G+e)=H(\G)$ is $2/2^K$. The remainder
of the proof is similar to the one for the K-SAT model.

We have established (\ref{eqconditionalsubadditivity}) and (\ref
{eqconditionalsubadditivityProb}). With this, the proof of
Proposition \ref{propsubadditivity} is complete.
\end{pf*}

Finally we give a simple proof of Corollary \ref{coroSATthreshold}.

\begin{pf*}{Proof of Corollary \ref{coroSATthreshold}}
Define $c^*_H=\sup\{c\ge0\dvtx  H(c)=c\}$. It suffices to show that
$H(c)<c$ for all $c>c^*_H$. For every $\delta>0$ we can
find $c_0\in(c,c+\delta)$ such that $H(c_0)<c_0$. By Lipshitz
continuity result of Theorem \ref{theoremMainResultER}
it follows that $H(c)\le H(c_0)+(c-c_0)<c$ for all $c>c_0$, and the
assertion is established.
\end{pf*}

%s5 #&#
\section{Proofs: Log-partition function in Erd{\"o}s--R\'{e}nyi
graphs}\label{sectionProofsLogPartition}
The following property serves as an analogue of (\ref
{eqobservation}). Given two hypergraphs $\G_i=([N],E_i)$, $i=1,2$, on
the same set of nodes $[N]$
for each one of the six models and each finite $\lambda$,
%
%e14 #&#
\begin{equation}
\label{eqlogZ-observation} \bigl|\log Z(\G_1)-\log Z(
\G_2)\bigr|=O\bigl(|E_1\Delta E_2|\bigr).
\end{equation}
This follows from the fact that adding (deleting) a hyperedge to (from)
a graph results in multiplying or dividing
the partition function by at most $\lambda$ for all models
except for the Ising and Independent set models. For the Ising model
the corresponding value is $\lambda^\beta$.
To obtain a similar estimate for the independent set, note that given a
graph $\G$ and an edge $e=(u,v)$ which is not in $\G$, we have
\[
Z(\G)=\sum_{e\subset I}\lambda^{|I|}+\sum
_{e\not\subset I}\lambda^{|I|},
\]
where in both sums we only sum over independent sets of $\G$. We claim that
\[
\sum_{e\subset I}\lambda^{|I|}\le\lambda\sum
_{e\not\subset
I}\lambda^{|I|}.
\]
Indeed, for every independent set in $\G$ containing $e=(u,v)$, delete
node $u$. We obtain a one-to-one mapping immediately leading
to the inequality. Finally, we obtain
\[
Z(\G)\ge Z(\G+e)=\sum_{e\not\subset I}\lambda^{|I|}
\ge{1\over
1+\lambda}Z(\G),
\]
where our claim was used in the second inequality.
Assertion (\ref{eqlogZ-observation}) then follows after taking logarithms.

The analogue of Theorem \ref{theoremsubadditivity} is the following result.
%
%th6 #&#
\begin{theorem}\label{theoremlogZ-subadditivity}
For every $1\le N_1, N_2\le N-1$ such that $N_1+N_2=N$ and every
$\lambda>1$,
\[
\E\bigl[\log Z\bigl(\G\bigl(N,\lfloor cN\rfloor\bigr)\bigr)\bigr]\ge\E
\bigl[\log
Z\bigl(\G(N_1,\rM_1)\bigr)\bigr]+\E\bigl[\log Z\bigl(
\G(N_2,\rM_2)\bigr)\bigr], %\label{eqlogZ-subadditive}
\]
where $\rM_1\distr\Bi(\lfloor cN\rfloor, N_1/N)$ and $\rM
_2\triangleq\lfloor cN\rfloor-\rM_1\distr\Bi(\lfloor cN\rfloor, N_1/N)$.
\end{theorem}

As before, we do not have independence of $\rM_j, j=1,2$.
Let us first show how this result implies Theorem
\ref{theoremlogZ-MainResultER}.
\begin{pf*}{Proof of Theorem \ref{theoremlogZ-MainResultER}}
Since $\rM_j$ have binomial distribution, using observation (\ref
{eqlogZ-observation}) and Theorem \ref{theoremlogZ-subadditivity},
we obtain
\begin{eqnarray*}%\label{eqsubadditive1LogZ}
&&\E\bigl[\log Z\bigl(\G\bigl(N,\lfloor cN\rfloor\bigr)\bigr)\bigr]
\\
&&\qquad\ge\E\bigl[\log Z\bigl(\G\bigl(N_1,\lfloor cN_1\rfloor
\bigr)\bigr)\bigr]+\E\bigl[Z\bigl(\G\bigl(N_2,\lfloor cN_2
\rfloor\bigr)\bigr)\bigr]-O(\sqrt{N}).
\end{eqnarray*}
Now we use Proposition \ref{propsubadditivityAppendix} in Appendix
\ref{appB}
for the case $\alpha=1/2$ to conclude that the
limit
\[
\lim_{N\rightarrow\infty} N^{-1}\E\bigl[\log Z\bigl(\G\bigl(N,
\lfloor cN\rfloor\bigr)\bigr)\bigr]\triangleq z(c)
\]
exists. Showing that this also implies the convergence of
$N^{-1}\E[\log Z(\G(N$, $\lfloor cN\rfloor))]$ to $z(c)$ w.h.p. again
is done using standard concentration results \cite{JansonBook} by
applying property (\ref{eqlogZ-observation}), and we skip the details.
The proof of continuity and monotonicity of $z(c)$ for relevant models
is similar to the one of $H(c)$.
\end{pf*}

Thus it remains to prove Theorem \ref{theoremlogZ-subadditivity}.
\begin{pf*}{Proof of Theorem \ref{theoremlogZ-subadditivity}}
We construct an interpolating sequence of graphs $\G(N,\lfloor
cN\rfloor,r), 0\le r\le\lfloor cN\rfloor$ exactly as in the previous
subsection.
We now establish the following analogue of Proposition
\ref{propsubadditivity}.
\end{pf*}

%pr3 #&#
\begin{prop}\label{proplogZ-subadditivity}
For every $r=1,\ldots,\lfloor cN\rfloor$,
%
%e15 #&#
\begin{equation}\label{eqlogZ-partialsubadditivity}
\E\bigl[\log Z\bigl(\G\bigl(N,\lfloor cN\rfloor,r\bigr)\bigr)\bigr] \ge
\E\bigl[
\log Z\bigl(\G\bigl(N,\lfloor cN\rfloor,r-1\bigr)\bigr)\bigr].
\end{equation}
\end{prop}
Let us first show how
Theorem \ref{theoremlogZ-subadditivity} follows from this proposition.
Observe that for disjoint union of two graphs $\G=\G_1+\G_2$, with
$\G=(V,E), \G_1=(V_1,E_1),\G_2=(V_2,E_2)$,
we always have $\log Z(\G)=\log Z(\G_1)+\log Z(\G_2)$. Theorem
\ref{theoremlogZ-subadditivity} then follows from Proposition
\ref{proplogZ-subadditivity}.

\begin{pf*}{Proof of Proposition \ref{proplogZ-subadditivity}}
Recall that $\G(N,\lfloor cN\rfloor,r-1)$ is obtained from $\G
(N,\lfloor cN\rfloor,r)$ by deleting a
hyperedge chosen u.a.r. independently from $r$ hyperedges $e_1,\ldots,e_r$
and adding a hyperedge $e$ either to nodes $[N_1]$ or to nodes $[N_2]$
with probabilities $N_1/N$ and $N_2/N$, respectively. Let as before $\G
_0$ be the graph obtained
after deleting but before adding a hyperedge, and let $Z_0=Z_0(\G_0)$
and $\mu_0=\mu_{0,\G_0}$ be the corresponding partition function and
the Gibbs measure, respectively.
In the case of K-SAT and NAE-K-SAT models we assume that $\G_0$
encodes the realizations of the random potentials as well.
We now show that conditional on any realization of the graph $\G_0$,
%
%e16 #&#
\begin{equation}\label{eqlogZ-conditionalsubadditivity}
\E\bigl[\log Z\bigl(\G\bigl(N,\lfloor cN\rfloor,r\bigr)\bigr)|\G_0
\bigr] \ge\E\bigl[\log Z\bigl(\G\bigl(N,\lfloor cN\rfloor,r-1\bigr
)\bigr)|
\G_0\bigr].
\end{equation}

The proof of (\ref{eqlogZ-conditionalsubadditivity}) is done on a
case-by-case basis, and it is very similar to the proof of
(\ref{eqconditionalsubadditivity}).\vspace*{9pt}

\textit{Independent sets}.
We have
\begin{eqnarray*}
&&
\E\bigl[\log Z\bigl(\G\bigl(N,\lfloor cN\rfloor,r\bigr)\bigr)|\G_0
\bigr]-\log Z_0 \\
&&\qquad= \E\biggl[\log\frac{Z(\G(N,\lfloor
cN\rfloor,r))}{Z_0}\Big|\G_0 \biggr]
\\
&&\qquad= \E\biggl[\log\frac{\sum_{I}\lambda^{|I|}-\sum_{I}1_{\{e\subset I\}
}\lambda^{|I|}}{\sum_{I}\lambda^{|I|}}\Big|\G_0\biggr]
\\
&&\qquad= \E\bigl[\log\bigl(1- \mu_0(e\subset I_0)\bigr)|
\G_0\bigr],
\end{eqnarray*}
where the sums $\sum_I$ are over independent sets only, and
$I_0$ denotes an independent set chosen randomly according to $\mu_0$.
Notice that since we are conditioning on graph $\G_0$, the only
randomness underlying the expectation operator is the randomness of the
hyperedge $e$ and the randomness of set $I_0$.
Note that $\mu_0(e\subset I_0)<1$ since $\mu_0(e\not\subset I_0)\geq
\mu_0(I_0=\varnothing)>0$.
Using the expansion $\log(1-x)=-\sum_{m\ge1}x^m/m$,
\begin{eqnarray*}
&&
\E\bigl[\log Z\bigl(\G\bigl(N,\lfloor cN\rfloor,r\bigr)\bigr)|\G_0
\bigr]-\log Z_0
\\
&&\qquad= - \E\Biggl[\sum_{k=1}^\infty
\frac{\mu_0(e\subset I_0)^k}{k}\Big|\G_0\Biggr]
\\
&&\qquad= -\sum_{k=1}^\infty\frac{1}{k} \E
\biggl[\sum_{I^1,\ldots,I^k}1_{\{
e\subset\bigcap_{j=1}^k I^j\}} \frac{\lambda^{\sum
_{j=1}^k|I^j|}}{Z_0^k}\Big|
\G_0\biggr]
\\
&&\qquad= -\sum_{k=1}^\infty\frac{1}{k}
\sum_{I^1,\ldots,I^k} \frac
{\lambda^{\sum_{j=1}^k|I^j|}}{Z_0^k} \E[ 1_{\{e\subset\bigcap_{j=1}^k
I^j\}} |
\G_0]
\\
&&\qquad= -\sum_{k=1}^\infty\frac{1}{k}
\sum_{I^1,\ldots,I^k} \frac
{\lambda^{\sum_{j=1}^k|I^j|}}{Z_0^k} \biggl(\frac{|\bigcap_{j=1}^k
I^j|}{N}
\biggr)^2,
\end{eqnarray*}
where the sum $\sum_{I^1,\ldots,I^k}$ is again over independent
subsets $I^1,\ldots,I^k$ of $\G_0$ only, and
in the last equality we have used the fact that $e$ is distributed u.a.r.
Similar calculation for $\log Z(\G(N,\lfloor cN\rfloor,r-1))$ that is
obtained by adding a hyperedge to nodes $[N_1]$ with probability
$N_1/N$, or to
nodes $[N_2]$ with probability $N_2/N$, gives
\begin{eqnarray*}
&&
\E\bigl[\log Z\bigl(\G\bigl(N,\lfloor cN\rfloor,r-1\bigr)\bigr)|
\G_0\bigr]-\log Z_0
\\
&&\qquad=-\sum_{k=1}^\infty\frac{1}{k}\sum
_{I^1,\ldots,I^k} \frac
{\lambda^{\sum_{j=1}^k|I^j|}}{Z_0^k} \biggl[\frac{N_1}{N}
\biggl(\frac{|\bigcap_{j=1}^k I^j\cap[N_1]|}{N_1}\biggr)^2\\
&&\qquad\quad\hspace*{114.5pt}{} +\frac
{N_2}{N}\biggl(
\frac{|\bigcap_{j=1}^k I^j\cap[N_2]|}{N_2}\biggr)^2 \biggr].
\end{eqnarray*}
Again using the convexity of $f(x)=x^2$ we obtain
\begin{eqnarray*}
&&\E\bigl[\log Z\bigl(\G\bigl(N,\lfloor cN\rfloor,r\bigr)\bigr)|\G_0
\bigr]-\log Z_0
\\
&&\qquad\geq\E\bigl[\log Z\bigl(\G\bigl(N,\lfloor cN\rfloor,r-1\bigr)\bigr)|
\G_0\bigr]-\log Z_0,
\end{eqnarray*}
and (\ref{eqlogZ-conditionalsubadditivity}) is established.\eject%\vspace*{9pt}

\textit{MAX-CUT}. Similarly to the independent set model, if
$\G(N,\lfloor cN\rfloor,r)$ is obtained from $\G_0$ by adding an edge
$(i,j)$ where $i,j$ are chosen uniformly at random, we have
\begin{eqnarray*}
&&
\E\bigl[\log Z\bigl(\G\bigl(N,\lfloor cN\rfloor,r\bigr)\bigr)|\G_0
\bigr]-\log Z_0
\\
&&\qquad= \E\biggl[\log\frac{Z(\G(N,\lfloor cN\rfloor,r))}{Z_0}\Big|\G_0\biggr]
\\
&&\qquad= \E\biggl[\log\frac{\sum_{x}1_{\{x_i=x_j\}}\lambda^{H(x)}+\lambda
\sum_{x}1_{\{x_i\neq x_j\}}\lambda^{H(x)}}{\sum_{x}\lambda
^{H(x)}}\Big|\G_0\biggr]
\\
&&\qquad= \log\lambda+ \E\biggl[\log\biggl(1- \biggl(1-\frac{1}{\lambda}\biggr)
\mu_0(x_i = x_j)\biggr)\Big|\G_0
\biggr].
\end{eqnarray*}

Since $\lambda>1$ we have $0<(1-\lambda^{-1})\mu_0(x_i = x_j)<1$
(this is where the condition $\lambda>1$ is used), implying
\begin{eqnarray*}
&&
\E\bigl[\log Z\bigl(\G\bigl(N,\lfloor cN\rfloor,r\bigr)\bigr)|\G_0
\bigr]-\log Z_0-\log\lambda
\\
&&\qquad= - \E\Biggl[\sum_{k=1}^\infty
\frac{(1-\lambda^{-1})^k\mu_0(x_i=
x_j)^k}{k}\Big|\G_0\Biggr]
\\
&&\qquad= -\sum_{k=1}^\infty\frac{(1-\lambda^{-1})^k}{k} \E
\biggl[\sum_{x_1,\ldots,x_k} \frac{\lambda^{\sum_{\ell=1}^kH(x_\ell
)}}{Z_0^k}1_{\{x^\ell_i=
x^\ell_j, \forall\ell\}}\Big|
\G_0\biggr]
\\
&&\qquad= -\sum_{k=1}^\infty\frac{(1-\lambda^{-1})^k}{k}\sum
_{x_1,\ldots,x_k} \frac{\lambda^{\sum_{\ell=1}^kH(x_\ell)}}{Z_0^k} \E
[1_{\{x^\ell
_i= x^\ell_j, \forall\ell\}} |
\G_0].
\end{eqnarray*}
Now for every sequence of vectors $x_1,\ldots,x_k$ introduce
equivalency classes on $[N]$. Given $i,k\in[N]$, say $i\sim k$ if
$x^\ell_i=x^\ell_k, \forall\ell=1,\ldots,k$. Namely, in every one
of the cuts defined by $x_\ell, \ell=1,\ldots,k$,
the nodes $i$ and $k$ belong to the same side of the cut.
Let $O_s, 1\le s\le J$ be the corresponding equivalency classes. For an
edge ${e=(i,j)}$ generated u.a.r., observe that
$\E[1_{\{x_i^\ell= x_j^\ell\forall\ell\}} |\G_0
]=\sum_{s=1}^J (\frac{|O_s|}{N} )^2$.
Thus
\begin{eqnarray*}
&&
\E\bigl[\log Z\bigl(\G\bigl(N,\lfloor cN\rfloor,r\bigr)\bigr)|\G_0
\bigr]-\log Z_0-\log\lambda
\\
&&\qquad
=-\sum_{k=1}^\infty\frac{(1-\lambda^{-1})^k}{k}\sum
_{x_1,\ldots,x_k} \frac{\lambda^{\sum_{\ell=1}^kH(\ell)}}{Z_0^k} \sum
_{s=1}^J \biggl(\frac{|O_s|}{N}
\biggr)^2
\end{eqnarray*}
and similarly,
\begin{eqnarray*}
&&
\E\bigl[\log Z\bigl(\G\bigl(N,\lfloor cN\rfloor,r-1\bigr)\bigr)|
\G_0\bigr]-\log Z_0-\log\lambda
\\
&&\qquad
= -\sum_{k=1}^\infty\frac{(1-{1}/{\lambda})^k}{k}\sum
_{x_1,\ldots,x_k} \frac{\lambda^{\sum_{\ell=1}^kH(\ell)}}{Z_0^k} \\
&&\qquad\quad{}\times\sum
_{s=1}^J \biggl(\frac{N_1}{N}\biggl(
\frac{|O_s\cap
[N_1]|}{N_1}\biggr)^2+\frac{N_2}{N}\biggl(
\frac{|O_s\cap[N_2]|}{N_2}\biggr)^2 \biggr).
\end{eqnarray*}
Using the convexity of the function $f(x)=x^2$,
we obtain (\ref{eqlogZ-conditionalsubadditivity}).\vspace*{9pt}

\textit{Ising}, \textit{coloring}, \textit{K-SAT and NAE-K-SAT}. The proofs of the
remaining cases are obtained similarly and are omitted. The condition
$\lambda>1$ is used
to assert positivity of $1-\lambda^{-1}$ in the logarithm expansion.
\end{pf*}

%s6 #&#
\section{Proofs: Random regular graphs}\label{sectionProofsReg}
For the proofs related to random regular graphs, we will need to work
with random ``nearly'' regular graphs. For this purpose,
given $N,r$ and $K$ such that $Nr/K$ is an integer and given any
positive integer $T\le Nr/K$,
let $\G(N,r,T)$ denote the graph obtained by creating
a size $T$ partial matching on $Nr$ nodes of the configuration model
uniformly at random and then projecting. For example,
if $T$ was $Nr/K$, then we would have obtained the random regular graph
$\G(N,r)$.

Our result leading to the proof of Theorem
\ref{theoremMainResultRegular} is as follows.
%
%th7 #&#
\begin{theorem}\label{theoremsubadditivityReg}
For every $N_1, N_2$ such that $N=N_1+N_2$ and $N_1r/K, N_2r/K$ are integers,
%
%e17 #&#
\begin{equation}
\label{eqsuperadOptimizationReg} \E\bigl[H\bigl(\G(N,r)\bigr)\bigr]\ge
\E\bigl[H
\bigl(\G(N_1,r)\bigr)\bigr]+\E\bigl[H\bigl(\G(N_2,r)\bigr)
\bigr]-O\bigl(N^{5/6}\bigr).
\end{equation}
\end{theorem}

\begin{pf}
Fix $N_1, N_2$ such that $N_1+N_2=N$ and
$N_1r/K$, $N_2r/K$ are integers.
Let us first prove Theorem \ref{theoremsubadditivityReg} for the
simpler case $\min_{j=1,2}N_j<40N^{5/6}$.
In this case starting from the graph\vspace*{1pt} $\G(N,r)$, we can obtain a
disjoint union of graphs $\G(N_j,r)$ via
at most $O(N^{5/6})$ hyperedge deletion and addition operations.
Indeed, suppose without the loss of generality that $N_1<40N^{5/6}$.
Delete all the hyperedges inside $[N_1]$ as well as all the hyperedges
connecting two parts. Then generate a random graph $\G(N_1,r)$
from scratch. Finally, complete a so-obtained partial matching in the
configuration model on $[N_2r]$ and project.
The total number of deleted and generated hyperedges is $O(N^{5/6})$,
and indeed we obtain a disjoint union
of graphs $\G(N_j,r), j=1,2$. Since the hyperedge deletion and
generation operation changes the value of $H$ by at most $O(N^{5/6})$,
then the proof of (\ref{eqsuperadOptimizationReg}) follows.

Now, throughout the remainder of the section we assume $\min
_{j=1,2}N_j\ge40N^{5/6}$.
Fix $T=Nr/K-\lfloor(1/K)N^{2/3}\rfloor$, and consider the graph
$\G(N,r,T)$.
Note that $Nr/K-T=O(N^{2/3})$.

We now describe
an interpolation procedure which interpolates between $\G(N,r,T)$ and
a union of certain two graphs on nodes
$[N_1]$ and $[N_2]$, each of which will be ``nearly'' regular. For every
integer partition $K=K_1+K_2$ such that $K_1,K_2\ge1$
let $T_{K_1,K_2}\le T$ be the (random) number of hyperedges which
connect parts $[N_1]$ and $[N_2]$ in $\G(N,r,T)$
and such that each connecting hyperedge has exactly $K_j$ nodes in part
$[N_jr]$ in the configuration model.
Let $T_0=\sum_{K_1,K_2\ge1: K_1+K_2=K}T_{K_1,K_2}$.
Observe that $T_0\le\min_{j=1,2}(N_jr)$.\vspace*{1pt}

Define $\G(N,T,0)=\G(N,r,T)$ and define $\G(N,T,t), 1\le t\le
T_{1,K-1}$, recursively as follows. Assuming $\G(N,T,t-1)$ is
already defined, consider the graph $\G_0$ obtained
from $\G(N,T,t-1)$ by deleting a hyperedge connecting $[N_1]$ and
$[N_2]$ chosen uniformly at random
from the collection of hyperedges which have exactly $1$ node in part
$[N_1r]$ and $K-1$ nodes in part $[N_2r]$
[from the remaining $T_{1,K-1}-(t-1)$ such hyperedges].
Then we construct
$\G(N,T,t)$ by adding a hyperedge to the resulting graph as follows:
with probability
$1/K$ a hyperedge is added to connect $K$ isolated nodes chosen
uniformly at random among the isolated nodes from the set $[N_1r]$.
With the remaining probability $(K-1)/K$
a hyperedge is added to connect $K$ isolated nodes chosen
uniformly at random among the isolated nodes from the set $[N_2r]$.
It is possible that at some point there are no $K$ isolated nodes
available in $[N_jr]$.
In this case we say that the interpolation procedure fails. In fact we
say that the interpolation procedure
fails if in either of the two parts the number of isolated nodes is
strictly less than $K$, even if the attempt
was made to add a hyperedge to a part where there is no shortage of
such nodes.

Thus we have defined an interpolation procedure
for $t\le T_{1,K-1}$. Assuming that the procedure did not fail for
$t\le T_{1,K-1}$, we
now define it for $T_{1,K-1}+1\le t\le T_{2,K-2}$ analogously: we
delete a randomly
chosen hyperedge connecting two parts such that the hyperedge has $2$
nodes in part $j=1$, and $K-2$ nodes in part $j=2$.
Then we add a hyperedge uniformly at random to part $j=1,2$ to connect
$K$ isolated nodes with probability $2/K$
and $(K-2)/K$, respectively. The failure of the interpolation is defined
similarly as above.
We continue this for all partitions $(K_1,K_2)$ until $(K-1,1)$, inclusive.
For the $(K_1,K_2)$ phase of the interpolation procedure the
probabilities are $K_1/K$ and $K_2/K$, respectively.

The interpolation procedure is particularly easy to understand in the
special case $K=2$. In this case $T_0=T_{1,1}$, and there is only one
phase in the interpolation procedure. In this phase every edge (which
is simply a pair of nodes) connecting sets $[N_1r]$ and $[N_2r]$ (if
any exists)
is deleted and replaces by an edge connecting two isolated nodes in
$[N_1r]$ with probability $1/2$, or two isolated nodes in $[N_2r]$ with
probability $1/2$ as well. One might note the difference of
probabilities $1/2$ and $1/2$ for the case of regular graphs vs.
$N_j/N, j=1,2$,
for the case of Erd{\"o}s--R\'{e}nyi graph. The difference stems from
the regularity assumption of the graph $\G(N,r)$.

Let $\mathcal{I}_t$ be the event that the interpolation succeeds for
the first $t$ steps, and let
$\mathcal{I}\triangleq\bigcap_{t\le T_0}\mathcal{I}_t$ denote the
event that the interpolation procedure succeeds for all steps.
For simplicity, even if the interpolation procedure fails in some step
$t'$, we still define $\G(N,T,t), t'\le t\le T_0$
to be the same graph as the first graph at which the interpolation
procedure fails, $\G(N,T,t)=\G(N,T,t')$.
It will be also convenient to define $\G(N,T,t)=\G(N,T,T_0)$ for
$T_0\le t\le\min_{j=1,2}(N_jr)$, whether the interpolation procedure
fails or not. This is done in order to avoid
dealing with graphs observed at a random ($T_0$) time, as opposed to
the deterministic time $\min_{j=1,2}(N_jr)$.

Provided that the interpolation procedure succeeds,
the graph
$\G(N,T,\break\min_{j=1,2}N_jr)$ is a disjoint union of two graphs on
$[N_j], j=1,2$, each ``close'' to being an $r$-regular random graph,
in some appropriate sense to be made precise later.

Our next goal is establishing the following analogue of
Proposition \ref{propsubadditivity}.
As in previous sections, let $\G_0$ denote the graph obtained from $\G
(N,T,t-1)$ after deleting a hyperedge connecting two parts, but before
a hyperedge is added to one of the parts, namely, before creating $\G(N,T,t)$,
conditioned on $t\le T_0$ and the event that the interpolation process
succeeds till $t - \bigcap_{t'\le t}\mathcal{I}_{t'}$.
If, on the other hand the interpolation procedure fails before $t$,
let $\G_0$ be the graph obtained at the last successful interpolation
step after the last hyperedge deletion.
Let $\Delta_i$ denotes the degree of the node $i\in[N]$ in the graph
$\G_0$, and let
\[
Z_j(t)\triangleq\sum_{i\in[N_j]}(r-
\Delta_{i})
\]
denote the number of isolated nodes in the $j$th part of the
configuration model for $\G_0$ for $j=1,2$.

%pr4 #&#
\begin{prop}\label{propsubadditivityReg}
For every $t\le\min_j N_jr$,
%
%e18 #&#
\begin{equation}\label{eqpartialsubadditivityReg}
\E\bigl[H\bigl(\G(N,T,t-1)\bigr)\bigr] \ge\E\bigl[H\bigl(\G(N,T,t)\bigr
)\bigr]-O
\biggl(\E\max_{j=1,2}{1\over Z_j(t)} \biggr).
\end{equation}
\end{prop}

\begin{pf}
The claim is trivial when $T_0+1\le t$, since the graph remains the same.
Notice also that
\[
\E\bigl[H\bigl(\G(N,T,t-1)\bigr)|\mathcal{I}_{t-1}^c
\bigr]=\E\bigl[H\bigl(\G(N,T,t)\bigr)|\mathcal{I}_{t-1}^c
\bigr],
\]
since the two graphs are identical, and thus the statement of the
proposition holds.

Now we will condition on the event $\mathcal{I}_t$.
We now establish a stronger result. Namely,
%
%e19 #&#
\begin{equation}\label{eqpartialsubadditivityRegCond}\qquad
\E\bigl[H\bigl(\G(N,T,t-1)\bigr)|\G_0\bigr]\ge\E\bigl[H\bigl(
\G(N,T,t)\bigr)|\G_0\bigr] -O \biggl(\max_{j=1,2}
{1\over Z_j(t)} \biggr).
\end{equation}
Observe that conditioned on obtaining graph $\G_0$, the graph $\G
(N,T,t-1)$ can be recovered from $\G_0$ in distributional sense
by adding a hyperedge connecting $K_1$ isolated nodes from $[N_1r]$ to
$K_2$ isolated nodes from $[N_2r]$, all chosen uniformly at random, and
then projecting.

We now conduct model-dependent, case-by-case analysis.\vspace*{9pt}

\textit{Independent sets.} In this case $K=2$, and the only
possibility is \mbox{$K_1=K_2=1$}.
As in the previous section, $O^*$ again denote the set of nodes in
$[N]$ which belong to \textit{every} largest independent set in $\G_0$. Then
in the case of creating graph $\G(N,T,t-1)$ from $\G_0$, the newly
added edge $e$ decreases $H$ by one if both
ends of $e$ belong to $O^*$, and leaves it the same otherwise.
The first event occurs with probability
\begin{eqnarray*}
&&
{\sum_{i_1\in O^*\cap[N_1],i_2\in O^*\cap[N_2]}(r-\Delta
_{i_1})(r-\Delta_{i_2})
\over
\sum_{i_1\in[N_1],i_2\in[N_2]}(r-\Delta_{i_1})(r-\Delta_{i_2})}
\\
&&\qquad={\sum_{i\in O^*\cap[N_1]}(r-\Delta_{i})
\over
\sum_{i\in[N_1]}(r-\Delta_{i})} {\sum_{i\in O^*\cap[N_2]}(r-\Delta_{i})
\over
\sum_{i\in[N_2]}(r-\Delta_{i})}.
\end{eqnarray*}
We now analyze the case of creating $\G(N,T,t)$. Conditioning on the
event that $e$ was added
to part $[N_jr]$, the value of $H$ decreases by one if and only if both
ends of $e$ fall into $O^*\cap[N_j]$. This occurs
with probability
\begin{eqnarray*}
%{\sum_{i\in O^*\cap[N_j]}(r-\Delta_{i_1})(r-\Delta_{i_2}-
&&
{ (\sum_{i\in[O^*\cap N_j]}(r-\Delta_i) )^2-\sum_{i\in
O^*\cap[N_j]}(r-\Delta_i)\over
(\sum_{i\in[N_j]}(r-\Delta_i) )^2-\sum_{i\in
[N_j]}(r-\Delta_i)}
\\
&&\qquad={ (\sum_{i\in[O^*\cap N_j]}(r-\Delta_i) )^2\over
(\sum_{i\in[N_j]}(r-\Delta_i) )^2}-O \biggl({1\over\sum_{i\in
[N_j]}(r-\Delta_i)} \biggr). %&={\Big(\sum_{i\in[O^*\cap N_j]}(r-
\end{eqnarray*}
%
%where in the last part we use the fact that initially there were $Z_j$
%free nodes in the part $[N_jr]$, and in each step the number goes down
%by one at most.
Therefore, the value of $H$ decreases by one with probability
\[
{1\over2}\sum_{j=1,2}
{ (\sum_{i\in[O^*\cap N_j]}(r-\Delta
_i) )^2\over
(\sum_{i\in[N_j]}(r-\Delta_i) )^2}-O \biggl(\max_{j=1,2}
{1\over Z_j(t)} \biggr)
\]
and stays the same with the remaining probability. Using the inequality
$(1/2)(x^2+y^2)\ge xy$, we obtain (\ref{eqpartialsubadditivityRegCond}).\vspace*{9pt}

\textit{MAX-CUT}, \textit{Ising}, \textit{coloring}. As in the proof of Theorem
\ref{theoremMainResultER},
we introduce equivalence classes $O_j^*\subset[N], 1\le j\le J$, on the
graph $\G_0$.
The rest of the proof is almost identical to the one for the
independent set model, and we skip the details. Notice that
in all of these cases we have $K=2$, and the interpolation phase has
only one stage corresponding to $(K_1,K_2)=(1,1)$.\vspace*{9pt}

\textit{K-SAT}. This is the first model for which $K>2$. We fix
$K_1,K_2\ge1$ such that $K_1+K_2=K$ and
further condition on the event that the graph $\G_0$ was created in
stage $(K_1,K_2)$.
As in the previous section, let $O^*$ be the set of frozen variables in
all optimal assignments of $\G_0$.
Reasoning as in the previous section, when we reconstruct graph $\G
(N,T,t-1)$ in the distributional sense by adding a
random hyperedge connecting $K_1$ nodes in $[N_1r]$ with $K_2$ nodes in
$[N_2r]$, the probability that the value
of $H$ remains the same (does not increase by one) is precisely
%
%e20 #&#
\begin{equation}
\label{eqleftH} {1\over2^K} \biggl[ {\sum_{i\in O^*\cap[N_1]}(r-\Delta_{i})
\over
\sum_{i\in[N_1]}(r-\Delta_{i})}
\biggr]^{K_1} \biggl[ {\sum_{i\in O^*\cap[N_2]}(r-\Delta_{i})
\over
\sum_{i\in[N_2]}(r-\Delta_{i})} \biggr]^{K_2}.
\end{equation}
Similarly, creating $\G(N,T,t)$ from $\G_0$ keeps the value of $H$
the same with probability
%
%e21 #&#
\begin{eqnarray}
\label{eqrightH}
&& {1\over2^K} {K_1\over K} \biggl[
{\sum_{i\in O^*\cap[N_1]}(r-\Delta_{i})
\over
\sum_{i\in[N_1]}(r-\Delta_{i})} \biggr]^{K} + {1\over2^K}
{K_2\over K} \biggl[ {\sum_{i\in O^*\cap[N_2]}(r-\Delta_{i})
\over
\sum_{i\in[N_2]}(r-\Delta_{i})} \biggr]^{K}
\nonumber\\[-8pt]\\[-8pt]
&&\qquad{}-O \biggl(\max_{j=1,2}{1\over\sum_{i\in[N_j]}(r-\Delta_i)} \biggr).
\nonumber
\end{eqnarray}
Applying Young's inequality, namely that $ab\le pa^{1/p}+qb^{1/q}$ for every $a,b\ge0$, $p+q=1, p,q>0$, with the choice
$p=K_1/K, q=K_2/K$,
\begin{eqnarray*}
a &=& \biggl[{\sum_{i\in O^*\cap[N_1]}(r-\Delta_{i})
\over
\sum_{i\in[N_1]}(r-\Delta_{i})} \biggr]^{K_1},
\\
b &=& \biggl[{\sum_{i\in O^*\cap[N_2]}(r-\Delta_{i})
\over
\sum_{i\in[N_2]}(r-\Delta_{i})} \biggr]^{K_2},
\end{eqnarray*}
and canceling $1/2^K$ on both sides, we obtain the result.\vspace*{9pt}

\textit{NAE-K-SAT}. The proof is similar to the one for K-SAT and for
NAE-K-SAT for the $\G(N,\lfloor cN\rfloor)$ model.
This completes the proof of the proposition.
\end{pf}

Our next step is to control the error term in (\ref{eqpartialsubadditivityReg}).

%le3 #&#
\begin{lemma}\label{lemmaerrorterm}
The interpolation procedure succeeds (event $\mathcal{I}$ holds) with
probability at least $1-O(N\exp(-N^\delta))$ for some $\delta>0$.
Additionally,
%
%e22 #&#
\begin{equation}
\label{eqerrorterm} \E\biggl[\sum_{1\le t\le T_0}\max
_{j=1,2}{1\over Z_j(t)} \biggr]=O\bigl(N^{2/5}
\bigr).
\end{equation}
\end{lemma}

\begin{pf}
Since $\G_0$ is obtained after deleting one hyperedge connecting two
parts, but before
adding a new hyperedge, then
$Z_j(t)\ge1$. A crude bound on the required expression is then $\E
[T_0]=O(\min N_j)$.
We have
$\E[Z_j(0)]=(N_j/N)N^{2/3}=N_j/N^{1/3}\ge40N^{1/2}$
since the initial number of isolated nodes was $Nr/K-T=N^{2/3}$
and $\min_j N_j\ge40N^{5/6}$.
Moreover, using a crude concentration bound $\pr
(Z_j(0)<(1/2)(N_j/N^{1/3})+K)=O(\exp(-N^{\delta_1}))$ for some
$\delta_1>0$. Observe that
$Z_j(t+1)-Z_j(t)=0$ with probability one if the interpolation procedure
failed for some $t'\le t$. Otherwise,
if $t$ corresponds to phase $(K_1,K_2)$, then
$Z_j(t+1)-Z_j(t)$ takes values $-K_j+K$ with probability $K_j/K$ and
$-K_j$ with the remaining probability.
This is because during the hyperedge deletion step, $Z_j(t)$
decreases by $K_j$, and during the hyperedge addition step, it
increases by $K$ or by zero with probabilities $K_j/K$ and $1-K_j/K$,
respectively.
In particular, $\E[Z_j(t+1)-Z_j(t)]=0$. The decision of whether to put
the hyperedge into part~$1$ or $2$ is made independently.
Since $t\le T_0\le N_j$, we conclude that for each $t\le T_0$ we have
$\pr(Z_j(0)-Z_j(t)>N_j^{3/5})=O(\exp(-N^{\delta_2}))$ for some
$\delta_2>0$. Here any choice of exponent strictly larger than $1/2$
applies, but for our purposes $3/5$ suffices.
It follows
that $Z_j(t)\ge(1/2)N_j/N^{1/3}+K-N_j^{3/5}$ for all $t$ with
probability $1-O(N_j\exp(-N^\delta))=1-O(N\exp(-N^\delta))$ for
$\delta=\min(\delta_1,\delta_2)$. The assumption $\min N_j\ge40
N^{5/6}$ implies that a weaker bound $\min N_j\ge32^{1/2}N^{5/6}$,
which translates into $(1/2)N_j/N^{1/3}-N_j^{3/5}\ge0$. Thus
$Z_j(t)\ge K$ for all $t$,
with probability $1-O(N\exp(-N^\delta))$, and therefore the
interpolation procedure succeeds.

Now ignoring term $K$ in the expression $(1/2)N_j/N^{1/3}+K-N_j^{3/5}$ and
using $T_0\le\min_j(N_jr)$, we obtain that with probability
$1-O(N\exp(-N^\delta))$, the expression inside the expectation
on the left-hand side of (\ref{eqerrorterm})
is at most
\[
{N_jr\over(1/2)N_jN^{-1/3}-N_j^{3/5}}={N_j^{2/5}r\over
(1/2) N_j^{2/5}N^{-1/3}-1}.
\]
The numerator is at most $N^{2/5}r$. Also the assumption $\min
N_j\ge40 N^{5/6}$ implies that the denominator
is at least $1$.
We conclude that the expression inside the expectation is at most
$N^{2/5}r$ with probability at least $1-O(N\exp(-N^\delta))$.
Since we also have $T_0\le Nr$ w.p.1, then
using a very crude estimate $O(N\exp(-N^\delta))=O(N^{-3/5})$,
and $NN^{-3/5}=N^{2/5}$,
we obtain the required result.
\end{pf}

As a corollary of Proposition \ref{propsubadditivityReg} and
Lemma \ref{lemmaerrorterm} we obtain
%
%co2 #&#
\begin{coro}\label{corointerpolationReg}
\[
\E\bigl[H\bigl(\G(N,T,0)\bigr)\bigr] \ge\E\Bigl[H\Bigl(\G\Bigl(N,T,\min
_{j=1,2}N_jr\Bigr)\Bigr)\Bigr]-O\bigl(N^{2/5}
\bigr).
\]
\end{coro}
Let us consider graph $\G(N,T,T_0)$. We further modify it by removing
all hyperedges which connect two parts $[N_j]$ of the graph,
if there are any such hyperedges left.
Notice that if the event $\mathcal{I}$ occurs, namely the
interpolation procedure succeeds, no further hyperedges need to be removed.
The resulting graph is a disjoint union of graphs obtained on nodes
$[N_1r]$ and $[N_2r]$ by adding a random size partial matching
uniformly at random. The actual size of these two matchings depends on
in the initial size of the partial matching within each part,
and also on
how many of $T_0$ hyperedges go into each part during the interpolation
steps, and how many
were removed in the final part (if any). We now obtain bounds on the
sizes of these matchings.

Recall $\min_j N_j\ge40N^{5/6}$.
We showed in the proof of Lemma \ref{lemmaerrorterm} that the
interpolation procedure succeeds with probability
$O(N\exp(-N^{\delta}))$ for some $\delta$. This coupled with the
fact that w.p.1, the number of hyperedges removed in the final stage is
at most $rN/K$,
gives us that the expected number of hyperedges removed in the final
stage is at most $O(N^2\exp(-N^{\delta}))$ which
(as a very crude estimate) is $O(N^{2/3})$.
Moreover, since the initial number of isolated nodes was $N^{2/3}$, and during the interpolation procedure the total
number of isolated nodes in the entire graph never increases, then the
total number of isolated nodes before the final removal of hyperedges
in $\G(N,T,T_0)$ is at most
$N^{2/3}$.
We conclude that the expected number of isolated nodes in the end of
the interpolation
procedure is $O(N^{2/3})$. Then we can complete uniform random
partial matchings on $[N_jr]$ to a perfect uniform random
matchings by adding at most that many hyperedges in expectation. The
objective value of $H$ changes by at most that much as well.
The same applies to $\G(N,r,T)$---we can complete the configuration
model corresponding to
this graph to a perfect matching on $Nr$ nodes by adding at most
$N^{2/3}$ hyperedges since $Nr/K-T=O(N^{2/3})$.
Coupled with Corollary \ref{corointerpolationReg} we then obtain
\[
\E\bigl[H\bigl(\G(N,r)\bigr)\bigr]\ge\E\bigl[H\bigl(\G(N_1,r)\bigr)
\bigr]+\E\bigl[H\bigl(\G(N_2,r)\bigr)\bigr]-O\bigl(N^{2/3}
\bigr)
\]
for the case $\min_j N_j\ge40N^{5/6}$.
This completes the proof of Theorem \ref{theoremsubadditivityReg}.
\end{pf}

\begin{pf*}{Proof of Theorem \ref{theoremMainResultRegular}}
The existence of the limit
\[
\lim_{N\rightarrow\infty, N\in r^{-1}K\Z_+} N^{-1}\E\bigl[H\bigl(\G
(N,r)\bigr)
\bigr]=H(r)
\]
follows immediately from Theorem \ref{theoremsubadditivityReg} and
Proposition \ref{propsubadditivityAppendix}
from Appendix \ref{appB}. Then the convergence w.h.p.
\[
\lim_{N\rightarrow\infty, N\in r^{-1}K\Z_+} N^{-1}H\bigl(\G(N,r)\bigr
)=H(r)
\]
follows once again using standard concentration results \cite{JansonBook}.
\end{pf*}

The proof of Theorem \ref{theoremlogZ-MainResultRegular} uses the
same interpolation as the one above, and the proof itself mimics the
one for
Theorem \ref{theoremlogZ-MainResultER}. For this reason, we omit the details.

%sA #&#
\begin{appendix}\label{app}

%sB #&#
\section{Proof of Lemma 1}%\protect\ref{lemmapN,m+1>CpN,m}
\label{appA}

We first assume K-SAT or NAE-K-SAT models.
Let us show for these models that there exists a constant $\omega\ge
1/2$ such that for \textit{every} graph and potential realization $(\G=(V,E),H)$
such that the problem is satisfiable [namely $H(\G)=|E|$], if a
randomly chosen hyperedge $e$ is added with a potential chosen
according to the model,
then
\[
\pr\bigl(H(\G+e)=|E|+1\bigr)\ge\omega.
\]
In other words, if the current graph is satisfiable, the new graph
obtained by adding a random hyperedge remains satisfiable with
probability at least $\omega$. Indeed, for example, for the case of
K-SAT, if the instance is satisfiable and $x$ is a satisfying assignment,
the added edge remains consistent with $x$ with probability at least
$\omega\triangleq1-1/2^K>1/2$. For the case of NAE-K-SAT it is
$\omega=1-1/2^{K-1}\ge1/2$.
We obtain that for every positive $M,m$ and recalling assumption
$\delta<1/2$,
\[
p(N,M+m)\ge\omega^m p(N,M)\ge\delta^m p(N,M),
\]
and the assertion is established.

The proof for the case of coloring is more involved. Given $0<\delta
<1/2$ we call a graph $\G$ on $N$ nodes $\delta$-unusual if it is
colorable, and in every coloring assignment there exists a color class
with size at least $(1-\delta)N$.
Namely, for every color assignment $x$ such that $H(x)=|E|$, there
exists $k\in[q^-]$ such that the cardinality of the set $\{i\in[N]\dvtx
x_i=k\}$ is at least $(1-\delta)N$.
We claim that
%
%eB.1 #&#
\begin{equation}
\label{equnusual} \pr\bigl(\G(N,M)\mbox{ is $\delta$-unusual}\bigr)\le(2
\delta)^M\exp\bigl(H(\delta)N+o(N)\bigr).
\end{equation}
The claim is shown using the first moment method---we will show that
the expected number of graphs with such a property is at most the
required quantity. Indeed,
given a subset $C\subset[N]$ such that $|C|\ge(1-\delta)N$, the probability
that the graph $\G(N,M)$ has proper coloring with all nodes in $C$
having the same color is at most $(1-(1-\delta)^2)^M<(2\delta)^M$,
since we must have that no edge falls within the class $C$. There are
at most
${N \choose\delta N}=\exp(H(\delta)N+o(N))$ choices for the subset
$C$. The claim then follows.

Now observe that if a graph $\G=(V,E)$ is colorable but not $\delta
$-unusual, then adding a random edge $e$, we obtain
$\pr(H(\G+e)=|E|+1)\ge\delta(1-\delta)$. Namely, in this case the
probability goes down by at most a constant factor.
We obtain
\begin{eqnarray*}
&&
p(N,M+1)\\
&&\qquad\ge\pr\bigl(H\bigl(\G(N,M+1)=M+1\bigr) |\G(N,M)
\mbox{ colorable, not
$\delta$-unusual} \bigr)
\\
&&\qquad\quad{}\times\pr\bigl(\G(N,M)\mbox{ colorable and not $\delta$-unusual}
\bigr)
\\
&&\qquad\ge\delta(1-\delta)\pr\bigl(\G(N,M)\mbox{ colorable and not $\delta$-unusual}
\bigr)
\\
&&\qquad\ge\delta(1-\delta)\pr\bigl(\G(N,M)\mbox{ colorable} \bigr)-\delta
(1-\delta)
\pr\bigl(\G(N,M)\mbox{ $\delta$-unusual} \bigr)
\\
&&\qquad\ge\delta(1-\delta)p(N,M)-\delta\pr\bigl(\G(N,M)\mbox{ $\delta$-unusual}
\bigr)
\\
&&\qquad\ge\delta(1-\delta)p(N,M)-\delta(2\delta)^M\exp\bigl(H(
\delta)N+o(N)\bigr),
\end{eqnarray*}
using the earlier established claim. Iterating this inequality, we
obtain for every $m\ge1$,
\begin{eqnarray*}
\hspace*{-4pt}&&p(N,M+m)
\\
\hspace*{-4pt}&&\qquad\ge\delta^m(1-\delta)^m p(N,M) -\delta(2
\delta)^M\exp\bigl(H(\delta)N+o(N)\bigr)\sum
_{0\le j\le m-1}\delta^m(1-\delta)^m
\\
\hspace*{-4pt}&&\qquad\ge\delta^m\pr\bigl(H\bigl(\G(N,M)=M\bigr) \bigr)-(2
\delta)^{M+1}\exp\bigl(H(\delta)N+o(N)\bigr),
\end{eqnarray*}
where $\sum_{0\le j\le m}\delta^m(1-\delta)^m\le1/(1-\delta)<2$ is
used in the last inequality. This completes the proof of the lemma.

%sC #&#
\section{Modified super-additivity theorem}\label{appB}
To keep the proof of our main results self-contained, we state and
prove the following proposition, used in proving several of the
theorems presented in the earlier sections. However,
B\'ela Bollob\'as and Zoltan F\"uredi kindly pointed out to us that the
following proposition is a special case of a more general and
% (which we had proved to keep our main contribution as well as the
%presentation self-contained),
classical theorem of de Bruijn and Erd\"os (see Theorem 22 on page 161
in \cite{deBE52}), which uses a weaker assumption on the additive term
in the near super-additivity hypothesis;
% Indeed, the following theorem more than suffices for our present
%purposes.
also see \cite{deBE51} and the Bollob\'as--Riordan percolation book
\cite{BR2006} for more recent applications of this useful tool.

%pr5 #&#
\begin{prop}\label{propsubadditivityAppendix}
Given $\alpha\in(0,1)$, suppose a nonnegative sequence $a_N$, \mbox{$N\ge1$}
satisfies
%
%eC.1 #&#
\begin{equation}
\label{eqsubadditivea1} a_N\ge a_{N_1}+a_{N_2}-O
\bigl(N^\alpha\bigr)
\end{equation}
for every $N_1,N_2$ s.t. $N=N_1+N_2$.
Then the limit $\lim_{N\to\infty}{a_N\over N}$ exists.
\end{prop}

\begin{pf}
It is convenient to define $a_N=a_{\lfloor N\rfloor}$ for every real,
but not necessarily integer value $N\ge1$.
It is then straightforward to check that property (\ref
{eqsubadditivea1}) holds when
extended to reals as well [thanks to the correction term $O(N^{\alpha})$].
Let
\[
a^*=\limsup_{N\rightarrow\infty}{a_N\over N}.
\]
Fix $\varepsilon>0$ and find $k$ such that $1/k<\varepsilon\le1/(k-1)$. Find
find $N_0=N_0(\varepsilon)$ such that $N_0^{-1}a_{N_0}\ge a^*-\varepsilon$,
$k^\alpha N_0^{\alpha-1}<\varepsilon$.
Clearly, such $N_0$ exists. Consider any $N\ge kN_0$.
Find $r$ such that $kN_02^r\le N\le kN_02^{r+1}$.
Applying (\ref{eqsubadditivea1}) iteratively with $N_1=N_2=N/2$ we obtain
\begin{eqnarray*}
a_N&\ge&2^r a_{N/2^r}-\sum
_{0\le l\le r-1}O \biggl(2^l\biggl({N\over
2^l}
\biggr)^\alpha\biggr)
\\
&=&2^ra_{N/2^r}-O \bigl(2^{(1-\alpha)r}N^\alpha
\bigr).
\end{eqnarray*}
Now let us find $i$ such that $(k+i)N_0\le N/2^r\le(k+i+1)N_0$. Note
$i\le k$. Again using (\ref{eqsubadditivea1})
successively with $N_0$ for $N_1$ and $N/2^r,
(N/2^r)-N_0,(N/2^r)-2N_0,\ldots$ for $N_2$, we obtain
\begin{eqnarray*}
a_{N/2^r}&\ge&(k+i)a_{N_0}-O \biggl(k\biggl(
{N\over2^r}\biggr)^{\alpha} \biggr)
\\
&\ge&(k+i)a_{N_0}-O \biggl(k\biggl({N\over2^r}
\biggr)^{\alpha} \biggr).
\end{eqnarray*}
Combining, we obtain
\begin{eqnarray*}
a_N&\ge&2^r(k+i)a_{N_0}-O
\bigl(2^{(1-\alpha)r}N^\alpha\bigr) -O \bigl(k2^{r(1-\alpha)}N^{\alpha}
\bigr)
\\
&=&2^r(k+i)a_{N_0} -O \bigl(k2^{r(1-\alpha)}N^{\alpha}
\bigr).
\end{eqnarray*}
Then
\begin{eqnarray*}
{a_N\over N}&\ge&{2^r(k+i)\over2^r(k+i+1)}{a_{N_0}\over N_0} -O
\bigl(k2^{r(1-\alpha)}N^{\alpha-1} \bigr)
\\
&\ge&\biggl(1-{1\over(k+i+1)}\biggr) \bigl(a^*-\varepsilon\bigr) -O
\bigl(k2^{r(1-\alpha)}N^{\alpha-1}\bigr)
\\
&\ge&(1-\varepsilon) \bigl(a^*-\varepsilon\bigr) -O \bigl(k2^{r(1-\alpha
)}N^{\alpha-1}
\bigr),
\end{eqnarray*}
where $1/k<\varepsilon$ is used in the last inequality. Now
\[
k2^{r(1-\alpha)}N^{\alpha-1}\le k2^{r(1-\alpha)}\bigl(k2^{r}N_0
\bigr)^{\alpha
-1}=k^\alpha N_0^{\alpha-1}<\varepsilon,
\]
again by the choice of $N_0$. We have obtained
\[
{a_N\over N}\ge(1-\varepsilon) \bigl(a^*-\varepsilon\bigr)-\varepsilon
\]
for all $N\ge N_0k$. Since $\varepsilon$ was arbitrary the proof is complete.
\end{pf}
\end{appendix}

\section*{Acknowledgments}

The authors are grateful for the insightful discussions with Silvio
Franz, Andrea Montanari, Lenka Zdeborov\'{a}, Florant Krzakala, Jeff
Kahn and James Martin. The authors thank Zoltan F\"uredi and B\'ela
Bollob\'as for bringing the deBruijn--Erd\"os theorem and other
relevant literature, mentioned in Appendix \ref{appB}, to the authors'
attention. The authors are especially thankful to anonymous referees
for helpful technical comments and notation suggestions. Authors also
thank Microsoft Research New England for the hospitality and the
inspiring atmosphere, in which this work began.

%suskaldyti doi

% imsref loaded by lrinkeviciute, 2013-04-24 12:26:20

\printaddresses


\begin{thebibliography}{29}
% BibTex style file: ims.bst, 2013-01-28
% Default style options (sort=0,type=number).
% Used options (sort=1,type=number).

%b1 #&#
\bibitem{AldousFavoriteProblemsNew}
\begin{bmisc}[auto:STB|2013/04/11|08:11:48]
\bauthor{\bsnm{Aldous},~\bfnm{D.}\binits{D.}}
\bhowpublished{Open problems. Preprint. Available at:
\texttt{\href{http://www.stat.berkeley.edu/\textasciitilde  aldous/Research/OP/sparse\_graph.html}{http://www.stat.berkeley.}
\href{http://www.stat.berkeley.edu/\textasciitilde  aldous/Research/OP/sparse\_graph.html}{edu/\textasciitilde  aldous/Research/OP/sparse\_graph.html}}}.
\bptok{imsref}%
\end{bmisc}
\endbibitem

%b2 #&#
\bibitem{AldousFavoriteProblems}
\begin{bmisc}[auto:STB|2013/04/11|08:11:48]
\bauthor{\bsnm{Aldous},~\bfnm{D.}\binits{D.}}
\bhowpublished{Some open problems. Preprint. Available at:
\texttt{\href{http://stat-www.berkeley.edu/users/aldous/Research/problems.ps}{http://stat-www.}
\href{http://stat-www.berkeley.edu/users/aldous/Research/problems.ps}{berkeley.edu/users/aldous/Research/problems.ps}}}.
\bptok{imsref}%
\end{bmisc}
\endbibitem

%b3 #&#
\bibitem{AldousSteelesurvey}
\begin{bincollection}[auto:STB|2013/04/11|08:11:48]
\bauthor{\bsnm{Aldous},~\bfnm{David}\binits{D.}} \AND
  \bauthor{\bsnm{Steele},~\bfnm{J.~Michael}\binits{J.~M.}}
(\byear{2004}).
\btitle{The objective method: Probabilistic combinatorial optimization and
  local weak convergence}.
In \bbooktitle{Probability on Discrete Structures}
(\beditor{\binits{H.}\bfnm{H.}~\bsnm{Kesten}}, ed.)
\bpages{1--72}.
\bpublisher{Springer}, \blocation{Berlin}.
\bptok{imsref}%
\end{bincollection}
\endbibitem

%b4 #&#
\bibitem{AlonSpencer}
\begin{bbook}[mr]
\bauthor{\bsnm{Alon},~\bfnm{Noga}\binits{N.}} \AND
  \bauthor{\bsnm{Spencer},~\bfnm{Joel~H.}\binits{J.~H.}}
(\byear{1992}).
\btitle{The Probabilistic Method}.
\bpublisher{Wiley}, \blocation{New York}.
\bid{mr={1140703}}
\bptok{imsref}%
\end{bbook}
\endbibitem

%b5 #&#
\bibitem{BollobasRandReg}
\begin{barticle}[mr]
\bauthor{\bsnm{Bollob{\'a}s},~\bfnm{B{\'e}la}\binits{B.}}
(\byear{1980}).
\btitle{A probabilistic proof of an asymptotic formula for the number of
  labelled regular graphs}.
\bjournal{European J. Combin.}
\bvolume{1}
\bpages{311--316}.
\bid{issn={0195-6698}, mr={0595929}}
\bptnote{check year}%
\bptok{imsref}%
\end{barticle}
\endbibitem

%b6 #&#
\bibitem{BollobasBook}
\begin{bbook}[mr]
\bauthor{\bsnm{Bollob{\'a}s},~\bfnm{B{\'e}la}\binits{B.}}
(\byear{2001}).
\btitle{Random Graphs}.
\bpublisher{Cambridge Univ. Press}, \blocation{Cambridge}.
\bptok{imsref}%
\end{bbook}
\endbibitem

%b7 #&#
\bibitem{BR2006}
\begin{bbook}[mr]
\bauthor{\bsnm{Bollob{\'a}s},~\bfnm{B{\'e}la}\binits{B.}} \AND
  \bauthor{\bsnm{Riordan},~\bfnm{Oliver}\binits{O.}}
(\byear{2006}).
\btitle{Percolation}.
\bpublisher{Cambridge Univ. Press}, \blocation{New York}.
\bid{mr={2283880}}
\bptok{imsref}%
\end{bbook}
\endbibitem

%b8 #&#
\bibitem{BollobasRiordanMetrics}
\begin{barticle}[mr]
\bauthor{\bsnm{Bollob{\'a}s},~\bfnm{B{\'e}la}\binits{B.}} \AND
  \bauthor{\bsnm{Riordan},~\bfnm{Oliver}\binits{O.}}
(\byear{2011}).
\btitle{Sparse graphs: Metrics and random models}.
\bjournal{Random Structures Algorithms}
\bvolume{39}
\bpages{1--38}.
\bid{doi={10.1002/rsa.20334}, issn={1042-9832}, mr={2839983}}
\bptok{imsref}%
\end{barticle}
\endbibitem

%b9 #&#
\bibitem{AminPC}
\begin{bmisc}[auto:STB|2013/04/11|08:11:48]
\bauthor{\bsnm{Coja-Oghlan},~\bfnm{A.}\binits{A.}}
(\byear{2012}).
\bhowpublished{Personal communication}.
\bptok{imsref}%
\end{bmisc}
\endbibitem

%b10 #&#
\bibitem{CopGamMohSor}
\begin{barticle}[mr]
\bauthor{\bsnm{Coppersmith},~\bfnm{Don}\binits{D.}},
  \bauthor{\bsnm{Gamarnik},~\bfnm{David}\binits{D.}},
  \bauthor{\bsnm{Hajiaghayi},~\bfnm{Mohammad~Taghi}\binits{M.~T.}} \AND
  \bauthor{\bsnm{Sorkin},~\bfnm{Gregory~B.}\binits{G.~B.}}
(\byear{2004}).
\btitle{Random {MAX} {SAT}, random {MAX} {CUT}, and their phase transitions}.
\bjournal{Random Structures Algorithms}
\bvolume{24}
\bpages{502--545}.
\bid{doi={10.1002/rsa.20015}, issn={1042-9832}, mr={2060633}}
\bptok{imsref}%
\end{barticle}
\endbibitem

%b11 #&#
\bibitem{deBE51}
\begin{barticle}[mr]
\bauthor{\bparticle{de} \bsnm{Bruijn},~\bfnm{N.~G.}\binits{N.~G.}} \AND
  \bauthor{\bsnm{Erd{\"o}s},~\bfnm{P.}\binits{P.}}
(\byear{1951}).
\btitle{Some linear and some quadratic recursion formulas. {I}}.
\bjournal{Indag. Math. (N.S.)}
\bvolume{13}
\bpages{374--382}.
\bid{mr={0047161}}
\bptok{imsref}%
\end{barticle}
\endbibitem

%b12 #&#
\bibitem{deBE52}
\begin{barticle}[mr]
\bauthor{\bparticle{de} \bsnm{Bruijn},~\bfnm{N.~G.}\binits{N.~G.}} \AND
  \bauthor{\bsnm{Erd{\"o}s},~\bfnm{P.}\binits{P.}}
(\byear{1952}).
\btitle{Some linear and some quadratic recursion formulas. {II}}.
\bjournal{Indag. Math. (N.S.)}
\bvolume{14}
\bpages{152--163}.
\bid{mr={0047162}}
\bptok{imsref}%
\end{barticle}
\endbibitem

%b13 #&#
\bibitem{DemboMontanariIsing}
\begin{barticle}[auto:STB|2013/04/11|08:11:48]
\bauthor{\bsnm{Dembo},~\bfnm{A.}\binits{A.}} \AND
  \bauthor{\bsnm{Montanari},~\bfnm{A.}\binits{A.}}
(\byear{2010}).
\btitle{Ising models on locally tree-like graphs}.
\bjournal{Ann. Appl. Probab.}
\bvolume{20}
\bpages{565--592}.
\bptok{imsref}%
\end{barticle}
\endbibitem

%b14 #&#
\bibitem{DemboMontanariSurvey}
\begin{barticle}[mr]
\bauthor{\bsnm{Dembo},~\bfnm{Amir}\binits{A.}} \AND
  \bauthor{\bsnm{Montanari},~\bfnm{Andrea}\binits{A.}}
(\byear{2010}).
\btitle{Gibbs measures and phase transitions on sparse random graphs}.
\bjournal{Braz. J. Probab. Stat.}
\bvolume{24}
\bpages{137--211}.
\bid{doi={10.1214/09-BJPS027}, issn={0103-0752}, mr={2643563}}
\bptnote{check related, check year}%
\bptok{imsref}%
\end{barticle}
\endbibitem

%b15 #&#
\bibitem{DemboMontanariSun}
\begin{bmisc}[auto:STB|2013/04/11|08:11:48]
\bauthor{\bsnm{Dembo},~\bfnm{A.}\binits{A.}},
  \bauthor{\bsnm{Montanari},~\bfnm{A.}\binits{A.}} \AND
  \bauthor{\bsnm{Sun},~\bfnm{N.}\binits{N.}}
(\byear{2011}).
\bhowpublished{Factor models on locally tree-like graphs.
Available at arXiv:\arxivurl{1110.4821}}.
\bptok{imsref}%
\end{bmisc}
\endbibitem

%b16 #&#
\bibitem{FranzLeone}
\begin{barticle}[mr]
\bauthor{\bsnm{Franz},~\bfnm{Silvio}\binits{S.}} \AND
  \bauthor{\bsnm{Leone},~\bfnm{Michele}\binits{M.}}
(\byear{2003}).
\btitle{Replica bounds for optimization problems and diluted spin systems}.
\bjournal{J. Stat. Phys.}
\bvolume{111}
\bpages{535--564}.
\bid{doi={10.1023/A:1022885828956}, issn={0022-4715}, mr={1972121}}
\bptnote{check related, check year}%
\bptok{imsref}%
\end{barticle}
\endbibitem

%b17 #&#
\bibitem{FranzLeoneToninelliRegular}
\begin{barticle}[auto:STB|2013/04/11|08:11:48]
\bauthor{\bsnm{Franz},~\bfnm{Silvio}\binits{S.}} \AND
  \bauthor{\bsnm{Leone},~\bfnm{Michele}\binits{M.}}
(\byear{2003}).
\btitle{Replica bounds for optimization problems and diluted spin systems}.
\bjournal{J. Phys. A Math. Gen.}
\bvolume{36}
\bpages{10967--10985}.
\bptok{imsref}%
\end{barticle}
\endbibitem

%b18 #&#
\bibitem{FranzMontanariPrivateCommunication}
\begin{bmisc}[auto:STB|2013/04/11|08:11:48]
\bauthor{\bsnm{Franz},~\bfnm{S.}\binits{S.}} \AND
  \bauthor{\bsnm{Montanari},~\bfnm{A.}\binits{A.}}
(\byear{2009}).
\bhowpublished{Personal communication}.
\bptok{imsref}%
\end{bmisc}
\endbibitem

%b19 #&#
\bibitem{Friedgut}
\begin{barticle}[mr]
\bauthor{\bsnm{Friedgut},~\bfnm{Ehud}\binits{E.}}
(\byear{1999}).
\btitle{Sharp thresholds of graph properties, and the {$k$}-sat problem}.
\bjournal{J. Amer. Math. Soc.}
\bvolume{12}
\bpages{1017--1054}.
\bnote{With an appendix by Jean Bourgain}.
\bid{doi={10.1090/S0894-0347-99-00305-7}, issn={0894-0347}, mr={1678031}}
\bptok{imsref}%
\end{barticle}
\endbibitem

%b20 #&#
\bibitem{gallagerLDPC}
\begin{bbook}[mr]
\bauthor{\bsnm{Gallager},~\bfnm{R.~G.}\binits{R.~G.}}
(\byear{1963}).
\btitle{Low-Density Parity-Check Codes}.
\bpublisher{MIT Press}, \blocation{Cambridge, MA}.
\bptok{imsref}%
\end{bbook}
\endbibitem

%b21 #&#
\bibitem{gamarnikLSAT}
\begin{barticle}[mr]
\bauthor{\bsnm{Gamarnik},~\bfnm{David}\binits{D.}}
(\byear{2004}).
\btitle{Linear phase transition in random linear constraint satisfaction
  problems}.
\bjournal{Probab. Theory Related Fields}
\bvolume{129}
\bpages{410--440}.
\bid{doi={10.1007/s00440-004-0345-z}, issn={0178-8051}, mr={2128240}}
\bptok{imsref}%
\end{barticle}
\endbibitem

%b22 #&#
\bibitem{gamarnikMaxWeightIndSet}
\begin{barticle}[mr]
\bauthor{\bsnm{Gamarnik},~\bfnm{David}\binits{D.}},
  \bauthor{\bsnm{Nowicki},~\bfnm{Tomasz}\binits{T.}} \AND
  \bauthor{\bsnm{Swirszcz},~\bfnm{Grzegorz}\binits{G.}}
(\byear{2006}).
\btitle{Maximum weight independent sets and matchings in sparse random graphs.
  {E}xact results using the local weak convergence method}.
\bjournal{Random Structures Algorithms}
\bvolume{28}
\bpages{76--106}.
\bid{doi={10.1002/rsa.20072}, issn={1042-9832}, mr={2187483}}
\bptok{imsref}%
\end{barticle}
\endbibitem

%b23 #&#
\bibitem{GuerraTon}
\begin{barticle}[mr]
\bauthor{\bsnm{Guerra},~\bfnm{Francesco}\binits{F.}} \AND
  \bauthor{\bsnm{Toninelli},~\bfnm{Fabio~Lucio}\binits{F.~L.}}
(\byear{2002}).
\btitle{The thermodynamic limit in mean field spin glass models}.
\bjournal{Comm. Math. Phys.}
\bvolume{230}
\bpages{71--79}.
\bid{doi={10.1007/s00220-002-0699-y}, issn={0010-3616}, mr={1930572}}
\bptnote{check related}%
\bptok{imsref}%
\end{barticle}
\endbibitem

%b24 #&#
\bibitem{JansonBook}
\begin{bbook}[mr]
\bauthor{\bsnm{Janson},~\bfnm{Svante}\binits{S.}},
  \bauthor{\bsnm{{\L}uczak},~\bfnm{Tomasz}\binits{T.}} \AND
  \bauthor{\bsnm{Rucinski},~\bfnm{Andrzej}\binits{A.}}
(\byear{2000}).
\btitle{Random Graphs}.
\bpublisher{Wiley}, \blocation{New York}.
\bid{doi={10.1002/9781118032718}, mr={1782847}}
\bptok{imsref}%
\end{bbook}
\endbibitem

%b25 #&#
\bibitem{JansonThomason}
\begin{barticle}[mr]
\bauthor{\bsnm{Janson},~\bfnm{Svante}\binits{S.}} \AND
  \bauthor{\bsnm{Thomason},~\bfnm{Andrew}\binits{A.}}
(\byear{2008}).
\btitle{Dismantling sparse random graphs}.
\bjournal{Combin. Probab. Comput.}
\bvolume{17}
\bpages{259--264}.
\bid{doi={10.1017/S0963548307008802}, issn={0963-5483}, mr={2396351}}
\bptok{imsref}%
\end{barticle}
\endbibitem

%b26 #&#
\bibitem{MontanariLDPCInterpolation}
\begin{barticle}[mr]
\bauthor{\bsnm{Montanari},~\bfnm{Andrea}\binits{A.}}
(\byear{2005}).
\btitle{Tight bounds for {LDPC} and {LDGM} codes under {MAP} decoding}.
\bjournal{IEEE Trans. Inform. Theory}
\bvolume{51}
\bpages{3221--3246}.
\bid{doi={10.1109/TIT.2005.853320}, issn={0018-9448}, mr={2239148}}
\bptok{imsref}%
\end{barticle}
\endbibitem

%b27 #&#
\bibitem{PanchenkoTalagrand}
\begin{barticle}[mr]
\bauthor{\bsnm{Panchenko},~\bfnm{Dmitry}\binits{D.}} \AND
  \bauthor{\bsnm{Talagrand},~\bfnm{Michel}\binits{M.}}
(\byear{2004}).
\btitle{Bounds for diluted mean-fields spin glass models}.
\bjournal{Probab. Theory Related Fields}
\bvolume{130}
\bpages{319--336}.
\bid{doi={10.1007/s00440-004-0342-2}, issn={0178-8051}, mr={2095932}}
\bptnote{check related}%
\bptok{imsref}%
\end{barticle}
\endbibitem

%b28 #&#
\bibitem{WormaldModelsRandomGraphs}
\begin{bincollection}[mr]
\bauthor{\bsnm{Wormald},~\bfnm{N.~C.}\binits{N.~C.}}
(\byear{1999}).
\btitle{Models of random regular graphs}.
In \bbooktitle{Surveys in Combinatorics, 1999 ({C}anterbury)}.
\bseries{London Mathematical Society Lecture Note Series}
\bvolume{267}
\bpages{239--298}.
\bpublisher{Cambridge Univ. Press}, \blocation{Cambridge}.
\bid{mr={1725006}}
\bptok{imsref}%
\end{bincollection}
\endbibitem

\end{thebibliography}
\end{document}